\documentclass{article}

\usepackage{amssymb,amsmath,amsthm,amscd,psfig}
\usepackage{mathrsfs}
\usepackage{graphicx}        

\def\C{\mathbb C}

\def\E{\mathbb E}
\def\H{\mathbb H}

\def\P{\mathbb P}
\def\Q{\mathbb Q}
\def\R{\mathbb R}
\def\S{\mathbb S}

\def\Y{\mathbb Y}
\def\Z{\mathbb Z}

\def\Gg{\gamma}

\def\GD{\Delta}
\def\GL{\Lambda}
\def\Gl{\lambda}
\def\Gm{\mu}
\def\Go{\omega}

\def\Gr{\varrho}

\def\GS{\Sigma}
\def\Gs{\sigma}

\def\Gt{\tau}
\def\Gz{\zeta}

\def\BE{\mathbf E}
\def\BF{\mathbf F}

\def\Bj{\mathbf j}
\def\BS{\mathbf S}
\def\BT{\mathbf T}
\def\BV{\mathbf V}

\def\BGg{\boldsymbol\gamma}

\def\BGm{\boldsymbol\mu}

\def\BGt{\boldsymbol\tau}
\def\BGz{\boldsymbol\zeta}

\def\kBGg{\ff{\boldsymbol\gamma}}

\def\kBGm{\ff{\boldsymbol\mu}}

\def\kBGt{\ff{\boldsymbol\tau}}

\def\cG{\mathcal G}

\def\cJ{\mathcal J}

\def\cP{\mathcal P}
\def\cR{\mathcal R}

\def\sE{\mathscr E}
\def\sF{\mathscr F}
\def\sG{\mathscr G}

\def\fM{\mathfrak M}

\def\fR{\mathfrak R}
\def\fS{\mathfrak S}

\def\dd{\partial}
\def\Bdd{\boldsymbol\dd}
\def\smin{\setminus}
\def\conj{conj}
\def\emp{\emptyset}
\def\eksi{\mathbf{id}}

\def\fix{\mathrm{Fix}}

\def\bfix{\mathbf{Fix}}
\def\bconj{\mathbf{Perm}}

\def\im{\sqrt{-1}}

\def\projc{\P^{1}(\C)}
\def\projr{\P^{1}(\R)}

\def\centerfig#1{\centerline{\psfig{figure=#1,silent=}}}

\def\rtree{\mathcal OTree(\Gs) }

\newcommand{\ve}[1]{\mathbf{#1}}

\newcommand{\real}[1]{#1(\R)}
\newcommand{\comp}[1]{#1(\C)}

\def\curve{(\GS;\ve{p})}

\newcommand{\cmod}[1]{\overline{M}_{#1}(\C)}
\newcommand{\cmodo}[1]{M_{#1}(\C)}
\newcommand{\umod}[1]{\overline{U}_{#1}(\C) }
\newcommand{\cdiv}[1]{D_{#1}(\C)}
\newcommand{\cdivc}[1]{\overline{D}_{#1}(\C)}

\newcommand{\rmod}[1]{\overline{M}^{\Gs}_{#1}(\R)}
\newcommand{\rmods}[1]{\overline{M}^{\Gs'}_{#1}(\R)}
\newcommand{\rmodo}[1]{M^{\Gs}_{#1}(\R)}

\newcommand{\csp}[1]{Conf_{(#1)}}
\newcommand{\cspp}[1]{{\widetilde{Conf}}_{(#1)}}

\newcommand{\cspq}[1]{\widetilde{C}_{(#1)}}
\newcommand{\csq}[1]{C_{(#1)}}
\newcommand{\csqc}[1]{\overline{C}_{#1}}

\newcommand{\dsq}[1]{B_{(#1)}}

\newcommand{\konj}[1]{\overline{#1}}
\newcommand{\hh}[1]{\hat{#1}}

\newcommand{\ff}[1]{\widetilde{#1}}

\newcommand{\cls}[1]{[\csqc{#1}]}

\newcommand{\opp}[1]{O_{#1}}

 \newtheorem{thm}{Theorem}[section]
 
 \newtheorem{lem}[thm]{Lemma}
 \newtheorem{prop}[thm]{Proposition}
 \theoremstyle{definition}
 \newtheorem{defn}[thm]{Definition}
 \theoremstyle{remark}
 \newtheorem{rem}[thm]{Remark}
 \newtheorem*{ex}{Example}
 \numberwithin{equation}{section}

\begin{document}
%
%
%
%
%
%
%
%
%
\title
{Graph homology of moduli space of pointed real curves of genus zero}

\author{\"Ozg\"ur Ceyhan}






\date{June 21, 2007}

\maketitle

\begin{abstract}
The moduli space $\rmod{\BS}$ parameterizes the isomorphism classes of
$\BS$-pointed stable real curves of genus zero which are invariant under
relabeling by the involution $\Gs$. This moduli space  is stratified
according to the degeneration types of  $\Gs$-invariant curves. The degeneration
types of $\Gs$-invariant curves are encoded by their dual trees with additional
decorations. We construct a combinatorial graph complex generated by the
fundamental classes of strata of $\rmod{\BS}$. We show that the homology  of
$\rmod{\BS}$ is isomorphic to the homology of our graph complex. We also give
a presentation of the fundamental group of  $\rmod{\BS}$.
\end{abstract}


\section{Introduction}

\subsection*{Moduli of pointed real curves of genus zero}
Let $\BS=\{s_1,\cdots,s_n\}$ be a finite set, and $\Gs$ be an involution
acting on it.  The moduli space $\rmod{\BS}$ of $\Gs$-invariant curves
is the fixed point set of the real structure
\begin{equation*}
c_{\Gs}: (\GS;p_{s_1},\cdots,p_{s_n}) \mapsto
(\konj{\GS};p_{\Gs(s_1)},\cdots,p_{\Gs(s_n)})
\end{equation*}
of the moduli space $\cmod{\BS}$ of $\BS$-pointed stable (complex)
curves of genus zero. The moduli space $\rmod{\BS}$ parameterizes
the isomorphism classes of  $\BS$-pointed stable real curves of genus zero
which are invariant under relabeling by the involution $\Gs$.

The moduli space of $\Gs$-invariant curves has recently attracted
attention in various contexts such as the representation theory of quantum
groups \cite{ehkr,hk,ka}, multiple $\zeta$-motives \cite{gm}, open
Gromov-Witten/Welschinger invariants  and related problems
\cite{c1,fu,liu,psw,s,w1,w2}. All of
these applications require descriptions of the (co)homology
groups and/or the fundamental group of $\rmod{\BS}$.

\subsection*{Graph complex of $\rmod{\BS}$}
The moduli space $\rmod{\BS}$ has a natural combinatorial stratification.
Each stratum $C_{\BGt}$ is determined by the degeneration type of
$\Gs$-invariant curves. Degeneration types of $\Gs$-invariant
curves are encoded  by trees $\BGt$ with corresponding decorations.

In this work, we introduce a combinatorial graph complex $\cG_\bullet$
where
\begin{eqnarray*}
\cG_d = \left( \bigoplus_{\BGt \mid \dim C_{\BGt} = d} \Z\ \cls{\BGt}\right)  / I_d.
\end{eqnarray*}
are Abelian groups generated by the relative fundamental classes $\cls{\BGt}$ of
the strata $\konj{C}_{\BGt}$ of $\rmod{\BS}$.  The additive relations between
the generators of $\cG_d$ are spanned by the relations introduced by Keel in
\cite{ke} (and studied in detail \cite{ge,km2}) and by some additional natural
relations which are known as the Cardy relations in physics literature. The
differential $\dd: \cG_d \to \cG_{d-1}$ is given by
\begin{eqnarray*}
\dd \ \cls{\BGt} = \sum_{\BGg /e = \BGt} \pm  \cls{\BGg}
\end{eqnarray*}
where the $\BGg$'s are the degeneration types of  $\Gs$-invariant curves
lying in the codimension one faces of $\csqc{\BGt}$. Then,  we prove that;

\smallskip
\noindent
{\bf Theorem.} $H_*(\rmod{\BS};\Z)$ is isomorphic to $H_*(\cG_\bullet)$.

\subsection*{Cohomology versus graph homology}
In their recent preprint \cite{ehkr}, Etingof and his collaborators calculated the cohomology
algebra $H^*(\rmod{\BS};\Q)$ in terms of generators and relations for the $\Gs = \mathrm{id}$
case.  Until this work, there was little known about the topology of
$\rmod{\BS}$ (see,  \cite{c,DJS,de,hk,ka}).

The graph homology, in a sense, treats the homology of the moduli space
$\rmod{\BS}$ in complementary directions to \cite{ehkr}: The graph complex
provides a recipe to calculate the homology of $\rmod{\BS}$ in  $\Z$ coefficients
for all possible  involutions $\Gs$, and it reduces to the cellular complex of the moduli
space $\rmod{\BS}$ for $\Gs=\mathrm{id}$ (see \cite{de,ka}).
In particular, when  $\Gs$ has no fixed elements, the moduli
space $\rmod{\BS}$ parameterizes $\Gs$-invariant curves having both non-empty
and empty real parts. To the best of our knowledge, this case has never been treated in the
literature, since most of the  recent literature on the moduli of real curves has been adopted
from the moduli space of pseudoholomorphic discs.

Moreover, our presentation is based on a stratification of $\rmod{\BS}$. All possible
homological relations between the strata of $\rmod{\BS}$ can be deduced from the
graph homology. One can easily  see that these relations include those
introduced in \cite{ge,km1,km2}
 (which are responsible for the higher associativity of $Comm_\infty$-algebras),
$A_\infty$-type relations (arising from the image of the differential of the graph complex)
and finally Cardy type relations. The graph homology of $\rmod{\BS}$  allows us
to define the quantum
cohomology of real  varieties as an open-closed homotopy algebra which
is very similar  to  \cite{ks}. Studying quantum cohomology of real  varieties
is one of the motivations of this work. This study
is presented in a subsequent paper \cite{c1}.

However, the  graph homology comes with a lack of product structure.
It is expected to give a product formula similar to the Kontsevich-Manin
description in \cite{km2}. 

\subsection*{Plan of this paper}
Section \ref{ch_n_curves} contains  a brief overview of basic facts about the
moduli space $\cmod{\BS}$ of $\BS$-pointed stable complex curves. In
Section \ref{ch_r_curves},  we define $\Gs$-invariant curves and provide a
combinatorial description of their  degeneration types.
In Section \ref{ch_r_moduli}, we introduce a set of
real structures on $\cmod{\BS}$, and their real parts $\rmod{\BS}$ as moduli
spaces of $\Gs$-invariant curves.  In Section \ref{sec_stratification}, we give the
stratification of $\rmod{\BS}$  according to degeneration types of
$\Gs$-invariant curves. In Section \ref{ch_rel_hom}, we calculate the
homology of the strata of $\rmod{\BS}$ relative to the union of their substrata
of codimension one and higher. In Section \ref{ch_graph},  we define the graph
complex and prove that its homology is isomorphic to the homology of
$\rmod{\BS}$. Finally, we give a presentation of the fundamental group of
$\rmod{\BS}$ by using the groupoid of paths transversal to codimension one strata.

\subsection*{Acknowledgements} I take this opportunity to express my deep
gratitude to my supervisors V. Kharlamov and M. Polyak for their
assistance and guidance in this work.

I also wish to thank to Yu.-I. Manin and M. Marcolli for their interest and
suggestions which have been invaluable for me.

I am thankful to  K. Aker, A. Bayer, O. Cornea, E. Ha, F. Lalonde,  A. Mellit,
D. Radnell and A. Wand for their comments and suggestions.


This paper was conceived during my stay  at  Israel Institute of
Technology-Technion, and its early draft was prepared at
Max-Planck-Institut f\"ur Mathematik (MPI), Bonn. Thanks are also due
to MPI, Technion, Institut de Recherhe Math\'{e}matique Avanc\'{e}e
de Strasbourg and Centre de Recherhes  Math\'ematiques for
their hospitality.

\subsection*{Notation/Convention}
We denote the finite set $\{s_1,\cdots,s_n\}$ by $\BS$, and the
symmetric group consisting of all permutations of $\BS$ by $\S_n$.
For an involution $\Gs \in \S_n$, we denote the subset $\{s \in \BS
\mid s=\Gs(s)\}$ by $\bfix(\Gs)$ and its complement $\BS \smin \bfix(\Gs)$
by  $\bconj(\Gs)$.

In this paper, the genus of all curves is zero except when the
contrary is stated. Therefore, we usually  omit mentioning the
genus of curves.


\section{Pointed complex curves and their moduli}
\label{ch_n_curves}

This section reviews the basic facts about pointed complex curves of
genus zero and their moduli spaces.

\subsection{Pointed  complex  curves}
\label{sec_n_curves}

An {\it $\BS$-pointed  curve} $\curve$ is a connected complex
algebraic curve $\GS$ with   distinct, smooth, {\it labeled points}
$\mathbf{p} = (p_{s_1},\cdots,p_{s_n}) \subset \GS$ satisfying the
following conditions:
\begin{itemize}
\item $\GS$ has only nodal singularities.%
\item The arithmetic genus of $\GS$ is equal to zero.
\end{itemize}
The nodal points and labeled points are called {\it special} points.

A {\it family of $\BS$-pointed  curves} over a complex
manifold $B(\C)$ is a proper, flat holomorphic map
$\pi_B: U_B(\C) \to B(\C)$ with
$n$ sections $p_{s_1},\cdots,p_{s_n}$ such that each geometric fiber
$(\GS(b);\mathbf{p}(b))$ is an $\BS$-pointed curve.

Two $\BS$-pointed curves  $\curve$ and $(\GS';\mathbf{p}')$ are {\it
isomorphic} if there exists a bi-holomorphic equivalence $\Phi: \GS
\to \GS'$ mapping $p_{s}$ to $p'_{s}$ for all $s \in \BS$.

An $\BS$-pointed curve is {\it stable} if its automorphism group is
trivial (i.e., on each irreducible component, the number of singular
points plus the number of labeled points is at least three).

\subsubsection{Graphs}

A {\it graph} $\Gg$ is a pair of finite sets of {\it vertices}
$\BV_{\Gg}$ and {\it flags} (or {\it half edges}) $\BF_{\Gg}$ with a
boundary map $\Bdd_{\Gg}: \BF_{\Gg} \to \BV_{\Gg}$ and an involution
$\Bj_{\Gg}: \BF_{\Gg} \to \BF_{\Gg}$ ($\Bj_{\Gg}^{2} =\eksi$). We
call $\BE_{\Gg}= \{(f_{1},f_{2}) \in \BF_{\Gg}^{2} \mid f_{1}=
\Bj_{\Gg}(f_{2})\ \& \ f_{1} \ne f_{2}\}$ the set of {\it edges},
and $\BT_{\Gg} = \{f \in \BF_{\Gg} \mid f = \Bj_{\Gg}(f)\}$ the set
of {\it tails}. For a vertex $v \in \BV_{\Gg}$, let $\BF_{\Gg}(v)
=\Bdd^{-1}_{\Gg}(v)$ and $|v| = |\BF_{\Gg}(v)|$ be the {\it valency}
of $v$.

We think of a graph $\Gg$ in terms of its {\it geometric
realization} $||\Gg||$ as follows: Consider the disjoint union of closed
intervals $\bigsqcup_{f_i \in \BF_{\Gg}} [0,1] \times f_i$, and
identify $(0, f_i)$ with $(0, f_j)$ if $\Bdd_\Gg(f_i) =
\Bdd_\Gg(f_j)$, and identify $(t, f_i)$ with $(1-t, \Bj_\Gg(f_i))$
for $t \in [0,1]$ and $f_i \ne \Bj_\Gg(f_i)$. The geometric realization of
$\Gg$ has a piecewise linear structure.

\begin{defn}
\label{def_tree} A {\it tree} is a graph whose geometric realization
is connected and simply-connected. If $|v|>2$ for all vertices, then
such a tree is called {\it stable}.
\end{defn}

There are only finitely many isomorphism classes of stable trees
whose set of tails $\BT_{\Gg} $ is equal to $\BS$ (see, \cite{m} or \cite{c}).
We call the isomorphism classes of such trees {\it $\BS$-trees}.

\subsubsection{Dual trees of $\BS$-pointed curves}

Let $\curve$ be an $\BS$-pointed stable curve and $\eta: \hat{\GS}
\to \GS$ be its normalization. Let
$(\hat{\GS}_v;\hat{\mathbf{p}}_v)$ be the following pointed stable
curve: $\hat{\GS}_{v}$ is a component of $\hat{\GS}$, and
$\hat{\mathbf{p}}_v$ is the set of points consisting of the
pre-images of special points on $\GS_{v} := \eta (\hat{\GS}_v)$. The
points $\hat{\ve{p}}_v = (p_{f_{1}},\cdots, p_{f_{|v|}})$ on
$\hat{\GS}_v$ are ordered by the elements $f_*$ in the set
$\{f_{1},\cdots, f_{|v|}\}$.

\begin{defn}
\label{def_dual_tree} The {\it dual tree} of an $\BS$-pointed stable
curve $\curve$ is an $\BS$-tree $\Gg$   where
\begin{itemize}
\item $\BV_{\Gg}$  is the set of components of $\hat{\GS}$.

\item $\BF_{\Gg}(v)$ is the set consisting of the pre-images of special
points in $\GS_v$.

\item $\Bdd_{\Gg}: f \mapsto v$ if and only if $p_f \in \hat{\GS}_v$.

\item $\Bj_{\Gg}: f \mapsto f$ if and only if $\eta(p_f)$ is a labeled
point, and $\Bj_{\Gg}: f_1 \mapsto f_2$ if and only if $p_{f_1} \in
\hat{\GS}_{v_1}$ and $p_{f_2} \in  \hat{\GS}_{v_2}$ are
the pre-images of a nodal point $\GS_{v_1} \cap \GS_{v_2}$.
\end{itemize}
\end{defn}

\subsubsection{Combinatorics of degenerations}
\label{sec_degen}

Let $\curve$ be an $\BS$-pointed stable curve whose  dual tree is
$\Gg$. Consider the deformations of a nodal point of  $\curve$. Such
a deformation of $\curve$ gives a {\it contraction} of an edge of
$\Gg$: Let $e=(f_{e},f^{e}) \in \BE_{\Gg}$ be the edge corresponding
to the nodal point which is deformed, and let  $\Bdd_{\Gg}(f_e) = v_{e},
\Bdd_{\Gg}(f^e) =v^{e}$.
Consider the equivalence relation $\thicksim$ on the set of
vertices, defined by: $v \thicksim v$ for all $v \in \BV_{\Gg} \smin
\{v_{e},v^{e}\}$, and  $v_{e} \thicksim v^{e}$. Then there is an
$\BS$-tree $\Gg/e$ whose set of  vertices is $\BV_{\Gg}/\thicksim$ and
whose set of flags is $\BF_{\Gg} \smin \{f_e,f^e\}$.  The boundary map and
involution of $\Gg/e$ are the restrictions of $\Bdd_{\Gg}$ and
$\Bj_{\Gg}$.

We use the notation $\Gg < \Gt$ in order to indicate that the
$\BS$-tree $\Gt$ is obtained by contracting a set of edges of $\Gg$.

\subsection{Moduli space of $\BS$-pointed curves}
\label{stratacomplex}

The moduli space $\cmod{\BS}$ is the space of  isomorphism classes
of $\BS$-pointed stable curves. This space is stratified according
to the degeneration types of $\BS$-pointed stable curves. The degeneration types
of $\BS$-pointed stable curves are combinatorially encoded by $\BS$-trees.
In particular, the principal stratum $\cmodo{\BS}$ parameterizes
$\BS$-pointed irreducible curves; i.e., it is associated to the one-vertex
$\BS$-tree. Therefore,  $\cmodo{\BS}$  is the quotient of the product
$(\projc)^{n}$ minus the diagonals $\GD = \bigcup_{k<l} \{(p_{s_1},
\cdots,p_{s_n})| p_{s_k} = p_{s_l} \}$ by $Aut(\projc) = PSL_2(\C)$.

\begin{thm}[Knudsen \& Keel, \cite{knu,ke}]

(a) For any $|\BS| \geq 3$,  $\cmod{\BS}$ is a smooth projective
algebraic variety of (real) dimension $2|\BS|-6$.

(b)  Any family of $\BS$-pointed stable curves over $B(\C)$ is induced
by a unique morphism $\kappa: B(\C) \to \cmod{\BS}$. The universal
family of curves $\umod{\BS}$ of $\cmod{\BS}$ is isomorphic to
$\cmod{\BS \cup \{s_{n+1}\}}$.

(c) For any $\BS$-tree $\Gg$, there exists a quasi-projective
subvariety  $\cdiv{\Gg} \subset \cmod{\BS}$ parameterizing the $\BS$-pointed
curves whose dual tree is given by $\Gg$. The subvariety $\cdiv{\Gg}$ is isomorphic to
$\prod_{v \in \BV_{\Gg}} \cmodo{\BF_\Gg(v)}$. Its (real) codimension
(in $\cmod{\BS}$) is $2|\BE_{\Gg}|$.

(d) $\cmod{\BS}$ is stratified by pairwise disjoint subvarieties
$\cdiv{\Gg}$. The closure $\cdivc{\Gg}$ of any stratum $\cdiv{\Gg}$  is stratified by
$\{\cdiv{\Gg'} \mid \Gg' \leq \Gg \}$.
\end{thm}

\subsection{Forgetful morphisms}
\label{sec_forgetful}

We say that $(\GS;p_{s_1},\cdots,p_{s_{n-1}})$ is obtained by
forgetting the labeled point $p_{s_n}$ of an $\BS$-pointed curve
$(\GS;p_{s_1},\cdots, p_{s_n})$. However, the resulting pointed
curve may well be unstable. This happens when the component
$\GS_{v}$ of $\GS$ supporting $p_{s_n}$ has only two additional
special points. In this case, we contract this component to its
intersection point(s) with the components adjacent to $\GS_{v}$.
With this {\it stabilization}, we extend this map to the whole space and
obtain $\pi_{\{s_n\}}: \cmod{\BS} \to \cmod{\BS'}$ where $\BS'= \BS
\smin \{s_n\}$.  There exists a canonical isomorphism $\cmod{\BS}
\to \umod{\BS'}$ commuting with the projections to $\cmod{\BS'}$. In
other words, $\pi_{\{s_n\}}: \cmod{\BS} \to \cmod{\BS'}$ can be
identified with the universal family of curves.

A very detailed study on the moduli space $\cmod{\BS}$ can be found
in Chapter 3.2 and 3.3 in \cite{m}, and also in \cite{ke,knu}.


\section{Pointed real curves of genus zero}
\label{ch_r_curves}

This section reviews some basic facts on $\BS$-pointed real curves
and their degeneration types.

\subsubsection{}

A {\it real structure} on a complex variety $X:=\comp{X}$ is an
anti-holomorphic involution  $c_X : X \to X$. The fixed point
set $\real{X} = \fix(c_X)$ of the involution $c_X$ is called the real
part of the variety $\comp{X}$  (or of the real structure $c_X$).

\subsection{$\Gs$-invariant curves and their families}
\label{sec_s_curve}

An $\BS$-pointed stable  curve $\curve$ is called {\it
$\Gs$-invariant}  if it admits a real structure $c_\GS: \GS \to \GS$
such that $c_\GS(p_s) =p_{\Gs(s)}$ for all $s \in \BS$.

Let $\pi_B: \comp{U_B} \to \comp{B}$ be a  family of $\BS$-pointed
stable curves with a pair of real structures
\begin{eqnarray*}
\begin{CD}
\comp{U_B}  @> {c_{U}} >>   \comp{U_B}      \\
@V{\pi_B}VV      @VV{\pi_B}V  \\
\comp{B}       @> {c_B}   >>    \comp{B}.
\end{CD}
\end{eqnarray*}
Such a family is called  {\it $\Gs$-equivariant} if the following conditions
are met;
\begin{itemize}
\item if $\pi^{-1}(b) = \GS$, then $ \pi^{-1}(c_B(b)) = \konj{\GS}$
for every $b \in B$;

\item $c_U:  z \in \GS = \pi^{-1}(b) \mapsto z \in \konj{\GS}= \pi^{-1}(c_B(b))$.
\end{itemize}
Here a complex curve $\GS$ is regarded as a pair
$\GS=(C,J)$, where $C$ is the underlying two-dimensional manifold, $J$ is a complex structure on $C$, and  $\overline{\GS} = (C,-J)$
is its complex conjugated pair.

\begin{rem}
If $\pi_B: \comp{U_B} \to \comp{B}$ is a $\Gs$-equivariant family, then
each $(\GS(b),\ve{p}(b))$ for $b \in \real{B}$ is $\Gs$-invariant. It follows from the fact
that the group of automorphisms of $\BS$-pointed stable curves is
trivial.
\end{rem}

\subsection{Combinatorial types of $\Gs$-invariant curves}
\label{sec_planar}

Degeneration types of $\Gs$-invariant curves are combinatorially
encoded by $\BS$-trees with additional decorations. This section
contains definitions of relevant combinatorial structures.

\subsubsection{Topological types of $\Gs$-invariant curves.}

Let $\curve$ be a $\Gs$-invariant curve. Each irreducible real
component $\GS_v$ of $\GS$ is isomorphic to $\projc$ with a real
structure that is either $z \mapsto \bar{z}$ or $z \mapsto
-1/\bar{z}$. Note that the real structure $z \mapsto -1/\bar{z}$ has empty
real part. Therefore, $\curve$ is either one of the following topological types:
\begin{itemize}
\item {\it Type 1}:  $\real{\GS}$ is a tree of real projective spaces having
only nodal singularities,

\item {\it Type 2}:  $\real{\GS}$ is the empty set,

\item {\it Type 3}:  $\real{\GS}$ is  a solitary nodal point.
\end{itemize}

This fact directly follows from the classification of real
structures on $\projc$ (see, for instance \cite{m2}) and Lefschetz's
fixed point theorem.

\begin{rem}
The terminology `type 1'  and `type 2' is commonly used in real algebraic
geometry. However, `type 3' has not yet been used in the literature. In this paper
it is used due to its convenience.
\end{rem}

\begin{rem}
If $\bfix(\Gs) \ne \emp$, then all $\Gs$-invariant curves are of
type 1. This follows from the fact that real parts of
$\Gs$-invariant curves cannot be the empty set (i.e., they can't be
type 2) and all special points must be distinct (i.e., they are not
type 3 either).
By contrast, $\Gs$-invariant curves can be of type
1, type 2 or type 3 when $\bfix(\Gs) = \emp$.
\end{rem}

\subsubsection{Oriented combinatorial types}

$\Gs$-invariant curves inherit additional structures on their sets of
special points. In this subsection, we introduce   `oriented' versions
of these structures for different topological types of $\Gs$-invariant
curves separately.

Let $\curve$ be a $\Gs$-invariant curve, and let   $\Gg$ be its dual tree.
We denote the set of real components $\{v \mid c_\GS (\GS_v) = \GS_v\}$
of $\curve$ by $\BV_\Gg^\R$.

\subsubsection*{$\Gs$-invariant curves of type 1.}

Let $(\hh{\GS};\hh{\ve{p}})$ be the normalization of a
$\Gs$-invariant curve $\curve$ of type 1. By identifying
$\hh{\GS}_v$ with  $\GS_v \subset \GS$,  we obtain a real structure
on $\hh{\GS}_v$ for a real component $\GS_v$. The real part
$\real{\hh{\GS}_v}$ of this real structure divides $\hh{\GS}_v$ into
two halves: two 2-dimensional open discs, $\GS^+_v$ and $\GS^-_v$,
having $\real{\hh{\GS}_v}$ as their common boundary in $\hh{\GS}_v$.
The real structure $c_\GS: \hh{\GS}_v \to \hh{\GS}_v$ interchanges
$\GS^\pm_v$, and  the complex orientations of $\GS^\pm_v$ induce two
opposite orientations on $\real{\hh{\GS}_v}$, called its {\it complex
orientations}.

If we fix a labeling of halves of $ \hh{\GS}_v$ by $\GS^\pm_v$ (or
equivalently,  if we orient   $\real{\hh{\GS}_v}$ with one of the
complex orientations),  then the set of pre-images of special
points $\hh{\ve{p}}_v \in  \hh{\GS}_v$
admits the following structures:

\begin{itemize}
\item {\it An oriented cyclic ordering on the set of special points
lying in $\real{\GS_v}$:} For any point $p_{f_i} \in (\hat{\ve{p}}_v
\cap \real{\hh{\GS}_v})$, there is a unique $p_{f_{i+1}} \in
(\hat{\ve{p}}_v \cap \real{\hh{\GS}_v})$ which follows the point
$p_{f_i}$ in the positive direction of $\real{\hh{\GS}_v}$ (the
direction which is determined by the complex orientation induced
by the orientation of $\GS^+_v$).

\item {\it An ordered two-partition of the set of special points
lying in $\GS_v \smin \real{\GS_v}$}. The subsets $\hat{\ve{p}}_v
\cap \GS_v^\pm$ of $\hat{\ve{p}}_v$ give an ordered  partition
of $\hat{\ve{p}}_v \cap (\GS_v \smin \real{\GS_v})$ into two disjoint
subsets.
\end{itemize}

The relative positions of the special points lying in $\real{\hh{\GS}_v}$
and the complex orientation of  $\real{\hh{\GS}_v}$  give an
{\it oriented cyclic ordering} on the corresponding labeling set
$\BF_\Gg^\R(v) :=
(\hat{\ve{p}}_v \cap \real{\hh{\GS}_v})$. We denote this oriented cyclic
ordering  by  $\{f_{r_1}\} < \cdots <
\{f_{r_{l-1}}\} < \{f_{r_{l}}\}$. Moreover, the
partition  $\{p_f \in \GS_v^\pm\}$ gives an {\it ordered two-partition}
$\BF_\Gg^\pm(v) := \{f \mid p_f \in \GS_v^\pm\}$ of
$\BF_\Gg (v) \smin \BF_\Gg^\R(v)$.

The {\it oriented combinatorial type} of a real component $\GS_v$ with a
fixed complex orientation is the following set of data:
\begin{eqnarray*}
o_v := \{ \mathrm{type}\ 1; \mathrm{two\ partition}\ \BF_\Gg^\pm(v);
\mathrm{oriented\ cyclic\ ordering\ on}\ \BF_\Gg^\R(v)\}.
\end{eqnarray*}

If we consider a $\Gs$-invariant curve $\curve$ with a  fixed complex
orientation at each real component, then the set of oriented combinatorial
types of real components
\begin{eqnarray*}
o := \{o_v \mid v \in \BV_\Gg^\R  \}
\end{eqnarray*}
is called an {\it oriented combinatorial type} of $\curve$.

\subsubsection*{$\Gs$-invariant curves of type 2.}

For a $\Gs$-invariant curve $\curve$ of type 2,  $\GS$ has
a unique real component, $\GS_v$, since the intersection of a pair of
real components of a $\Gs$-invariant curve must be
a real point, and $\GS$ has none. Moreover, $\bfix(\Gs)$ must be the empty set since
$\real{\GS}=\emp$.

In this case, the {\it oriented combinatorial type} of $\curve$ is the
following set of data:
\begin{eqnarray*}
o := \{ \mathrm{type}\ 2; \BV_\Gg^\R = \{v\} \}.
\end{eqnarray*}

\subsubsection*{$\Gs$-invariant curves of type 3.}

For a $\Gs$-invariant curve $\curve$ of type 3,  the real part
$\real{\GS}$ of $\curve$ divides $\GS$ into two connected pointed
complex curves $(\GS^\pm;\ve{p}^\pm)$ having $\real{\GS}$ as their
intersection point. We denote the sets of components of
$(\GS^\pm;\ve{p}^\pm)$ by $\BV_\Gg^\pm$, and the sets of their flags
$\bigcup_{v \in\BV_\Gg^\pm} \Bdd_\Gg^{-1} (v)$ by $\BF_\Gg^\pm$.

An {\it oriented combinatorial type} of $\curve$ is the following
set of data:
\begin{eqnarray*}
o := \{ \mathrm{type}\ 3; \mathrm{two\ partitions}\ \BV_\Gg^\pm \
\mathrm{and}\ \BF_\Gg^\pm\}.
\end{eqnarray*}

\subsubsection{Unoriented combinatorial types}

The definition of oriented combinatorial types requires additional
choices (such as complex orientations) which are not determined by real structures of
$\Gs$-invariant curves.  By identifying the oriented combinatorial types for
such different choices, we obtain {\it unoriented combinatorial types}
of $\Gs$-invariant curves.

\subsubsection*{$\Gs$-invariant curves of type 1.}

For each real component $\GS_v$ of a $\Gs$-invariant curve $\curve$
of type 1, there are two possible ways of choosing $\GS_v^+$ in
$\hh{\GS}_v$. These two different choices give the {\it opposite}
oriented combinatorial types $o_v$ and $\bar{o}_v$. Namely, the
oriented combinatorial type $\bar{o}_v$ is obtained from $o_v$ by
reversing the cyclic ordering of $\BF_\Gg^\R(v)$ and swapping
$\BF_\Gg^+(v)$ and $\BF_\Gg^-(v)$.

The {\it unoriented combinatorial type} of a real
component $\GS_v$ of $\curve$  is the pair of opposite oriented
combinatorial types $u_v := \{o_v,\bar{o}_v\}$. The set of unoriented
combinatorial  types of real components
\begin{eqnarray*}
u := \{u_v \mid v \in \BV_\Gg^\R  \}
\end{eqnarray*}
is called the {\it unoriented combinatorial type} of $\curve$ .

\subsubsection*{$\Gs$-invariant curves of type 2.}

For a $\Gs$-invariant curve $\curve$ of type 2, the {\it unoriented
combinatorial type} is the same set of data as the oriented
combinatorial type i.e., $u := o = \{\text{type\ 2};
\BV_\Gg^\R = \{v\}\}$.

\subsubsection*{$\Gs$-invariant curves of type 3.}

For a $\Gs$-invariant curve $\curve$ of type 3, there are two
possible ways of choosing $(\GS^+;\ve{p}^+)$ in $\curve$. These
 choices give two {\it opposite} oriented combinatorial
types, $o$ and $\bar{o}$. Namely, the oriented combinatorial type
$\bar{o}$ is obtained from $o$ by swapping $\BV_\Gg^+$ and
$\BV_\Gg^-$, and swapping $\BF_\Gg^+$ and $\BF_\Gg^-$.

The {\it unoriented combinatorial type} of $\curve$ is the pair of
opposite oriented combinatorial types $u := \{o,\bar{o}\}$.

\subsubsection{Dual trees of $\Gs$-invariant curves}

The combinatorial types of $\Gs$-invariant curves can be encoded on
their dual trees.

\subsubsection*{O-planar trees.}
Let $\curve$ be an $\Gs$-invariant curve and let $\Gg$ be its dual
tree.

An {\it oriented locally planar (o-planar) structure} on $\Gg$
is a set of data which encodes an oriented combinatorial type of
$\curve$. O-planar structures for different
topological types are explicitly given as follows:
\begin{itemize}
\item For a $\Gs$-invariant curve of type 1,

\begin{itemize}
\item $\curve$ is of type 1 (i.e., $\real{\GS}$ is a tree of real projective spaces).

\item $\BV_{\Gg}^\R \subset \BV_\Gg$ is the set of real components
of $\GS$ (i.e., the set of {\it real vertices}).

\item $\BF_{\Gg}^\R (v) \subset \BF_\Gg(v)$ is the set of the
pre-images of special points in $\real{\GS_v}$. (i.e., the set of
{\it real flags} adjacent to the real vertex $v \in \BV_{\Gg}^\R$).

\item $\BF_{\Gg}^\R(v)$ carries an oriented cyclic ordering  for every
$v \in \BV_{\Gg}^\R$.

\item $\BF_\Gg(v) \smin \BF_{\Gg}^\R(v)$ carries an ordered two-partition
for every $v \in \BV_{\Gg}^\R$.
\end{itemize}

\item For a $\Gs$-invariant curve of type 2,

\begin{itemize}
\item $\curve$ is of type 2 (i.e., $\real{\GS}=\emp$).

\item $\BV_{\Gg}^\R = \{v\} \subset \BV_\Gg$ is the set of
real components of $\GS$ (i.e., the set of {\it real vertices}
contains only one element).
\end{itemize}

\item For a $\Gs$-invariant curve of type 3,

\begin{itemize}
\item $\curve$ is of type 3 (i.e., $\real{\GS}$ is a point).

\item $e=(f_e,f^e)$ is the edge corresponding to the solitary nodal
point of $\GS$.

\item  $\BF_\Gg$ and $\BV_\Gg$ carry two partitions
$\BF_\Gg^\pm$ and $\BV_\Gg^\pm$   respectively.
\end{itemize}
\end{itemize}

We denote $\BS$-trees $\Gg,\Gt,\Gm$ with o-planar structures
by $(\Gg,o),(\Gt,o),(\Gm,o)$ or
by bold Greek characters with tilde
$\kBGg,\kBGt,\kBGm$. When it is necessary to indicate different
o-planar structures on the same $\BS$-tree, we use indices in parentheses (e.g., $\kBGt_{(i)}$).

\subsubsection*{Notations.}

For each vertex  $v \in \BV_{\Gg}^\R$ (resp. $v \in \BV_{\Gg} \smin
\BV_{\Gg}^\R$) of an o-planar tree $\kBGg$, we associate a subtree
$\kBGg_{v}$ (resp. $\Gg_v$) which is given by $\BV_{\Gg_{v}} =
\{v\}, \BF_{\Gg_{v}} = \BF_{\Gg}(v)$, $\Bj_{\Gg_{v}} = \eksi$,
$\Bdd_{\Gg_{v}} = \Bdd_\Gg$, and by the  o-planar structure $o_v$ of
$\kBGg$ at the vertex $v \in \BV_\Gg^\R$.

A pair of vertices $v,\bar{v} \in \BV_\Gg \smin \BV_\Gg^\R$ is said to
be {\it conjugate} if $c_\GS(\GS_v) = \GS_{\bar{v}}$. Similarly, we
call a pair of flags $f,\bar{f} \in \BF_\Gg \smin \BF_\Gg^\R$ {\it
conjugate} if $c_\GS$ swaps the corresponding special points.

To each o-planar tree $\kBGg$ of type 1,
we associate the subsets of vertices $\BV_\Gg^\pm$ and flags
$\BF_\Gg^\pm$  as follows:
Let $v_1 \in \BV_\Gg \smin \BV_\Gg^\R$, and let $v_2 \in \BV_\Gg^\R$
be the closest vertex to $v_1$ in $||\Gg||$. Let $f(v_1) \in
\BF_\Gg(v_2)$ be in the shortest path  connecting the vertices $v_1$
and $v_2$. The sets $\BV_\Gg^\pm$ are the subsets of vertices $v_1 \in
\BV_\Gg \smin \BV_\Gg^\R$ such that the corresponding flags $f(v_1)$ are
respectively in $\BF_\Gg^\pm(v_2)$. The subsets of
flags $\BF_\Gg^\pm$ are defined as
$\Bdd_\Gg^{-1}(\BV_\Gg^\pm)$.

\subsubsection*{U-planar trees.}
A u-planar structure on the dual tree $\Gg$ of $\curve$ is the set of
data encoding the unoriented combinatorial type of $\curve$. It is
given by

\begin{eqnarray*}
u := \left\{
\begin{array}{ll}
\{(\Gg_v,o_v), (\Gg_v,\bar{o}_v) \mid v \in \BV_\Gg^\R\} &
\mathrm{if}\ \curve \ \mathrm{is\ of\ type\ 1}, \\
\{ \mathrm{type}\ 2; \BV_\Gg^\R = \{v\} \}  & \mathrm{if}\ \curve \
\mathrm{is\ of\ type\ 2}, \\
\{\mathrm{special\ real\ edge\ e=(f_e,f^e)}\} & \mathrm{if}\ \curve\
\mathrm{is\ of\ type\ 3}.
\end{array}
\right.
\end{eqnarray*}

We denote $\BS$-trees  $\Gg,\Gt,\Gm$ with u-planar structures by
$(\Gg,u),(\Gt,u),(\Gm,u)$ or simply by  bold Greek characters
$\BGg,\BGt,\BGm$. O-planar planar trees $\kBGg,\kBGt,\kBGm$
give representatives of u-planar trees $\BGg,\BGt,\BGm$ respectively.

\subsubsection{Contraction morphism of o/u-planar trees}

\subsubsection*{Contraction morphism of o-planar trees}
Consider a family of $\Gs$-invariant curves
which is a deformation of a real node of the central fiber
$(\GS(b_0),\ve{p}(b_0))$ with a given oriented combinatorial type.
Let $\kBGt,\kBGg$ be the o-planar trees associated respectively to a
generic fiber $(\GS(b),\ve{p}(b))$ and  the central fiber
$(\GS(b_0),\ve{p}(b_0))$ of this family. Let $e$ be the edge corresponding
to the nodal point that is deformed. We say that $\kBGt$ is
obtained by {\it contracting} the edge $e$ of  $\kBGg$, and to
indicate this we use the notation  $\kBGg < \kBGt$.

\subsubsection*{Contraction morphism of u-planar trees}
\begin{figure}[htb]
\centerfig{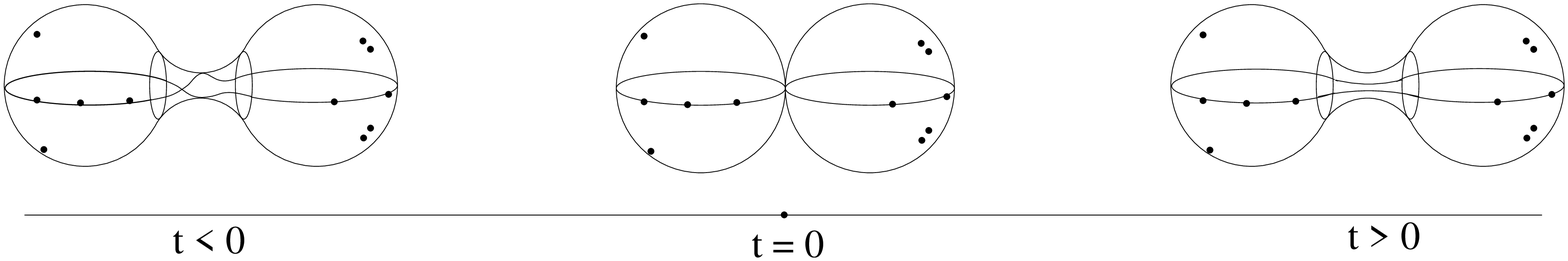,height=.8in} %
\caption{Two possible deformations of a real nodal point.}%
\label{deformation}
\end{figure}

The definition of contraction morphisms of u-planar trees is
the same as that of contraction morphisms of o-planar trees. By contrast,
the contraction of an edge of an u-planar tree is not a well-defined
operation: We can think of a deformation of a real node as
the family $\{x \cdot y = t \mid t \in \R \}$. According to the sign
of the deformation parameter $t$, we obtain two different
unoriented combinatorial types of $\Gs$-invariant curves,
see Figure \ref{deformation}. Different u-planar trees
$\BGg_{(i)}$ that are obtained  by contraction of the same edge of
 $\BGt$ correspond to different signs of deformation parameters.

Further details of contraction morphisms for o/u-planar trees can be
found in \cite{c}.


\section{Moduli of $\Gs$-invariant curves}
\label{ch_r_moduli}

The moduli space $\cmod{\BS}$ comes equipped with a natural real
structure. The involution
\begin{equation}
\label{eqn_p_r_str}
c: \curve  \mapsto(\overline{\GS};\ve{p})
\end{equation}
 gives the principal real structure of $\cmod{\BS}$.

On the other hand, the permutation group $\S_n$ acts on $\cmod{\BS}$
via relabeling: For each $\Gr \in \S_n$, there is an holomorphic map
$\psi_{\Gr}$ defined by $\curve \mapsto (\GS;\Gr(\ve{p})) := (\GS;
p_{\Gr(s_1)}, \cdots,p_{\Gr(s_n)})$. For each involution $\Gs \in
\S_n$, we have an additional real structure
\begin{equation}
\label{eqn_r_str}
c_{\Gs} := c \circ \psi_{\Gs}: \curve \mapsto
(\overline{\GS};\Gs(\ve{p}))
\end{equation}
of $\cmod{\BS}$. The real part $\rmod{\BS}$ of the real structure
$c_{\Gs}: \cmod{\BS} \to \cmod{\BS}$ gives the moduli space of
$\Gs$-invariant curves:

\begin{thm}[Ceyhan, \cite{c}]
\label{thm_r_moduli}

(a) For any $|\BS| \geq 3$, $\rmod{\BS}$ is a smooth projective real
manifold of dimension $|\BS|-3$.

(b) Any  $\Gs$-equivariant family $\pi_B: \comp{U_B} \to \comp{B}$
of  $\BS$-pointed stable curves is induced by a unique pair of real
morphisms
\begin{eqnarray}
\begin{CD}
\comp{U_B}  @> \hat{\kappa} >>   \umod{\BS} \\
@V{\pi_S}VV      @VV{\pi}V   \\
\comp{B}       @> \kappa >>       \cmod{\BS}.          \nonumber
\end{CD}
\end{eqnarray}

(c) Let $\fM_\Gs(\C)$ be the contravariant functor that sends each real
variety $(\comp{B},c_B)$ to the set of $\Gs$-equivariant families of
curves over $B$. The moduli functor $\fM_\Gs(\C)$ is represented by the
real variety $(\cmod{\BS},c_\Gs)$.

(d) Let $\real{\fM_\Gs}$ be the contravariant functor that sends each
real analytic manifold $R$ to the set of families of $\Gs$-invariant
curves over $R$. The moduli functor $\real{\fM_\Gs}$ is represented by
the real part $\rmod{\BS}$ of $(\cmod{\BS},c_\Gs)$.
\end{thm}

\begin{rem}
The group of holomorphic automorphisms of $\cmod{\BS}$ that respect
its stratification is isomorphic to $\S_n$. Therefore, the real structures
preserving the stratification of $\cmod{\BS}$  are of  the form (\ref{eqn_r_str}) (see \cite{c}).

However, we don't know whether there exist real structures other than
(\ref{eqn_r_str}), since the whole group of holomorphic automorphisms
$Aut(\cmod{\BS})$ is not necessarily isomorphic to $\S_n$. For example,
the automorphism group of $\cmod{\BS}$ is $PSL_2(\C)$ when $|\BS|=4$.

It is believed that $Aut(\cmod{\BS}) \cong \S_n$ for $|\BS| \geq 5$. In fact,
it is true for $|\BS|=5$ and a proof can be found in \cite{dol}.
To the best of our knowledge, there is no systematic exposition of
$Aut(\cmod{\BS})$ for $|\BS| >5$.
\end{rem}


\section{Stratification of moduli  of $\Gs$-invariant curves}
\label{sec_stratification} The  moduli space $\rmod{\BS}$ can be
stratified according to degeneration  types of $\Gs$-invariant
curves. This stratification is given in terms of  spaces of
$\Z_2$-equivariant point configurations on $\projc$ and moduli
spaces of complex curves $\cmod{\BF}$.

\subsection{Spaces of $\Z_2$-equivariant point configurations on $\projc$ }
\label{sec_conf_sp}

Let $z := [z:1]$ be  an affine coordinate on $\projc$. Consider the
upper half-plane $\H^+ = \{z \in \projc \mid \Im{(z)} >0 \}$ (resp.
lower half plane $\H^- = \{z \in \projc \mid \Im{(z)} <0 \}$) as a
half of the $\projc$ with respect to $z \mapsto \bar{z}$, and the
real part $\projr$ as its boundary. We denote  the compactified
disc  $\H^{+} \cup \projr$ by $\H$.

Let $\BF$ be a finite set and $\Gr$ be an involution acting on
$\BF$. Denote the fixed point set of $\Gr$ by $\BF^\R$.

\subsubsection{$\Z_2$-equivariant configurations on $(\projc, z \mapsto \konj{z})$}

The configuration space of $k =|\BF \smin \BF^\R|/2$ distinct pairs
of conjugate points in $\H^+ \bigsqcup \H^-$ and $l=|\BF^\R|$
distinct points in $\projr$ is
\begin{eqnarray*}
\cspp{\BF,\Gr} &:=& \{ (z_{f_{1}},\cdots,z_{f_{2k}};x_{g_1},\cdots,x_{g_l}) \mid
                        z_f \in \H^+ \bigsqcup \H^- \ \mathrm{for}\ f \in \BF \smin \BF^\R,\\
               &  &     z_{f} =      z_{f'} \Leftrightarrow f = f',
                        z_{f} = \bar{z}_{f'} \Leftrightarrow f' = \Gr(f) \  \& \\
               &  &     x_g \in \projr \ \mathrm{for}\ g \in \BF^\R, \ x_g =  x_{g'}
                        \Leftrightarrow g=g' \}.
\end{eqnarray*}
The number of connected components of $\cspp{\BF,\Gr}$ is
$2^{k} (l-1)!$  when $l \geq 2$, and $2^{k} $ when $l =0,1$.
 They are all pairwise diffeomorphic; the natural diffeomorphisms are given
 by $\Gr$-invariant relabelings.

The action of  $SL_2(\R)$ on $\H$ is given by
\begin{eqnarray}
SL_2(\R) \times \H \to \H,\ \ (\GL,z) \mapsto \GL(z)
= \frac{az+b}{cz+d},\ \
\GL =  \left(
\begin{array}{cc}
a   &   b \\
c   &   d \\
\end{array}
\right) \in SL_2(\R). \nonumber
\end{eqnarray}
 It induces an isomorphism $SL_2(\R)/ \pm I
\to Aut(\H)$. The  automorphism group $Aut(\H)$ acts on
$\cspp{\BF,\Gr}$ by
\begin{equation*}
\GL: (z_{f_1},\cdots,z_{f_{2k}};x_{g_1},\cdots,x_{g_l}) \mapsto
     (\GL(z_{f_1}),\cdots,\GL(z_{f_{2k}}); \GL(x_{g_1}),\cdots,\GL(x_{g_l})).
\end{equation*}
This action preserves each connected component of $\cspp{\BF,\Gr}$.
It is free when $2k+l \geq 3$, and it commutes with
diffeomorphisms given by $\Gr$-invariant relabelings. Therefore,
the quotient space $\cspq{\BF,\Gr} :=\cspp{\BF,\Gr}  / Aut(\H)$ is a
manifold of dimension $2k+l-3$ whose connected components are
pairwise diffeomorphic.

In addition to the automorphisms considered above, there is a
diffeomorphism of $\cspp{\BF,\Gr}$ which is given in affine
coordinates as follows
\begin{equation*}
e: (z_{f_1},\cdots,z_{f_{2k}}; x_{g_1},\cdots,x_{g_l})  \mapsto
   (-z_{f_1},\cdots,-z_{f_{2k}};-x_{g_1},\cdots,-x_{g_l}).
\end{equation*}
Consider the quotient space $\csp{\BF,\Gr} := \cspp{\BF,\Gr}/(e)$.
The diffeomorphism $e$ commutes  with each $\Gr$-invariant
relabeling and the normalizing action of $Aut(\H)$. Therefore, the
quotient space $\csq{\BF,\Gr} := \csp{\BF,\Gr} /Aut(\H)$ is a
manifold of dimension $2k+l-3$, its connected components are
diffeomorphic to the components of $\cspq{\BF,\Gr}$, and, moreover,
the quotient map $\cspq{\BF,\Gr} \to \csq{\BF,\Gr}$ is a trivial
double covering.

\subsubsection{$\Z_2$-equivariant configurations on $(\projc, z \mapsto -1/\konj{z})$}

Let $\BF^\R =\emp$. We denote by $Conf_{(\BF,\Gr)}^\emp$ the space
of configurations of $|\BF|$ labeled points on $\projc$ that are
invariant under the real structure $z \mapsto -1/\bar{z}$:
\begin{equation*}
Conf_{(\BF,\Gr)}^\emp :=\{(z_{f_1},\cdots,z_{f_{2k}}) \mid z_{\Gr(f)} = -1/\bar{z}_{f}\}.
\end{equation*}
The group of automorphisms of $\projc$ commuting with the real structure $z \mapsto -1/\bar{z}$ is
\begin{eqnarray}
Aut(\projc,\conj) \cong SU_2 := \left\lbrace  \left(
\begin{array}{cc}
a   &   b \\
-\bar{b}   &   \bar{a} \\
\end{array}
\right) \in SL_2(\C) \right\rbrace.  \nonumber
\end{eqnarray}
The group $Aut(\projc,\conj)$ acts naturally  on
$Conf_{(\BF,\Gr)}^\emp$. For $|\BF| \geq 4$, the action is free and
the quotient $\dsq{\BF,\Gr} := Conf_{(\BF,\Gr)}^\emp /
Aut(\projc,\conj)$ is a $2k-3$ dimensional connected manifold.

\subsubsection{Connected Components of $\rmodo{S}$}

Set $(\BF,\Gr) = (\BS,\Gs)$. Each connected component of
$\csq{\BS,\Gs}$ for $\bfix(\Gs) \ne \emp$ (resp. $\csq{\BS,\Gs}
\cup \dsq{\BS,\Gs}$ for $\bfix(\Gs)=\emp$) is associated to an
unoriented combinatorial type of $\Gs$-invariant curves, and
each unoriented combinatorial type is  given by a one-vertex u-planar
tree $\BGg$. We denote the connected components of $\csq{\BS,\Gs}$
(and $\csq{\BS,\Gs} \cup \dsq{\BS,\Gs}$) by  $C_{\BGg}$.

Every $\Z_2$-equivariant point configuration defines a
$\Gs$-invariant curve. Hence, we define
\begin{equation} \label{eqn_diffeo}
\Xi: \bigsqcup_{\BGg : |\BV_{\Gg}|=1} C_{\BGg} \to \rmodo{\BS}
\end{equation}
which maps  $\Z_2$-equivariant point configurations to the
corresponding isomorphisms classes of irreducible $\Gs$-invariant
curves.

\begin{lem}
\label{lem_conf_sp}

(a) The map $\Xi$ is a diffeomorphism.

(b) Let $|\bconj(\Gs)|=2k $ and $\bfix(\Gs)=l$. The configuration
space $C_{\BGg}$ is diffeomorphic to
\begin{itemize}
\item $((\H^+)^k \smin \GD) \times \R^{l-3}$ when $l >2$,

\item $((\H^+ \smin \{\im\})^{k-1}  \smin \GD )\times \R^{l-1}$
when $l=1,2$,

\item $((\H^+ \smin \{\im, \im/2 \})^{k-2}  \smin \GD )\times \R$
when $l=0$ and $\BGg$ of type 1,

\item $((\projc \smin \{-1,-1/2,2,1\})^{k-2}  \smin (\GD\cup\GD^c))
\times \R$ when $l=0$ and  $\BGg$ is of type 2.
\end{itemize}
Here, $\GD$ is the  union of all diagonals where $z_s = z_{s'}\ (s \ne
s')$, and $\GD^c$ is the union of all cross-diagonals where $z_s =
-\frac{1}{\bar z_{s'}} \ (s \ne s')$.
\end{lem}

The proof of Lemma \ref{lem_conf_sp} can be found in \cite{c}.

\subsection{Configuration spaces of o/u-planar trees and stratification of $\rmod{\BS}$}
\label{sec_ou_p}

We associate a product of configuration spaces  $C_{\kBGg_v}$ and
moduli spaces of pointed complex curves $\cmod{\BF_\Gg(v)}$ to each
o-planar tree $\kBGg$:
\begin{eqnarray*}
C_{\kBGg} =
     \left\lbrace
      \begin{array}{ll}
          \prod_{v \in \BV_{\Gg}^\R} C_{\kBGg_v} \times
             \prod_{\{v,\bar{v}\} \in \BV_\Gg \smin \BV_\Gg^\R} \cmodo{\BF_\Gg(v)}
             &  \mathrm{if}\ \BGg \ \mathrm{is\ of  \ type \ 1}, \\
          C_{\BGg_{v_r}} \times
             \prod_{\{v,\bar{v}\} \subset \BV_\Gg \smin \BV_{\Gg}^\R} \cmodo{\BF_\Gg(v)}
             &  \mathrm{if}\ \BGg \ \mathrm{is\ of  \ type \ 2},  \\
          \prod_{\{v,\bar{v}\} \subset \BV_\Gg \smin \BV_{\Gg}^\R} \cmodo{\BF_\Gg(v)}
             &  \mathrm{if}\ \BGg \ \mathrm{is\ of  \ type \ 3}.
      \end{array}
     \right.  \\
\end{eqnarray*}
Here, the products $\prod \cmodo{\BF_\Gg(v)}$ run over unordered
conjugate pairs of vertices.

For each u-planar $\BGg$, we first choose an o-planar representative
$\kBGg$, and then we set $C_{\BGg} := C_{\kBGg}$. Note that $C_{\BGg}$
does not depend on the o-planar representative.

\begin{thm}[Ceyhan, \cite{c}]
\label{thm_strata}

(a) $\rmod{\BS}$ is stratified by $C_{\BGg}$.

(b) A stratum $C_{\BGg}$ is contained in the boundary of
$\csqc{\BGt}$ if and only if  $\BGt$ is obtained by contracting an
invariant set of edges of $\BGg$. The codimension of $C_{\BGg}$ in
$\csqc{\BGt}$ is $|\BE_{\Gg}|-|\BE_{\Gt}|$.
\end{thm}

\begin{ex}
\label{exa_strata} The first nontrivial example is $\cmod{\BS}$
where $\BS=\{s_1,s_2,s_3,s_4\}$. There are three real structures:
$c_{\Gs_1}$, $c_{\Gs_2}$, $c_{\Gs_3}$, where
\begin{eqnarray*}
\Gs_1 = \mathrm{id}, \
\Gs_2 = \left(
          \begin{array}{cccc}
             s_1 & s_2 & s_3 & s_4 \\
             s_2 & s_1 & s_3 & s_4 \\
          \end{array}
         \right) \
\mathrm{and} \
\Gs_3 = \left(
          \begin{array}{cccc}
             s_1 & s_2 & s_3 & s_4 \\
             s_3 & s_4 & s_1 & s_2 \\
          \end{array}
         \right).
\end{eqnarray*}
Note that, all other real structures which preserve the stratification of
$\cmod{\BS}=\projc$ are conjugate to $c_{\Gs_i}, i=1,2,3$.

The main stratum  $\real{M_{\BS}^{\Gs_1}}$ of
$\real{\konj{M}_{\BS}^{\Gs_1}}$ is the configuration space
of four distinct points in $\projr$ up to the action of $PSL_{2}(\R)$.
Due to Lemma \ref{lem_conf_sp}, $\Gs_1$-invariant curves
$\curve \in \real{M_{\BS}^{\Gs_1}}$ are identified with tuples
 $(0,x_{s_2},1,\infty)$ where  $x_{s_2} \in \projr \smin \{0,1,\infty\}$.
Hence, the main stratum $\real{M_{\BS}^{\Gs_1}}$ is $\projr \smin \{0,1,\infty\}$,
and its compactification $\real{\konj{M}_{\BS}^{\Gs_1}}$ is $\projr$.  Three
intervals of $\real{M_{\BS}^{\Gs_1}}$ are the three configuration
spaces $C_{\BGt_{(i)}}$ and three  points  are the configuration
spaces  $C_{\BGg_{i}}$. The o-planar representatives of $\BGt_{(i)}$
and $\BGg_i$ are given in Fig. \ref{fig_rmoduli1}.

\begin{figure}[htb]
\centerfig{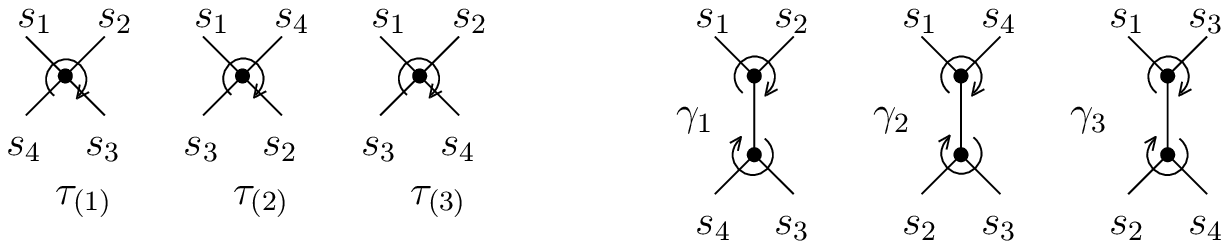,height=.9in} %
\caption{O-planar representatives of u-planar trees of the strata
of $\real{\konj{M}_{\BS}^{\Gs_1}}$.} %
\label{fig_rmoduli1}
\end{figure}

The main stratum $\real{M_{\BS}^{\Gs_2}}$
of $\real{\konj{M}_{\BS}^{\Gs_2}}$ is the space of distinct
configurations of two points in $\projr$ and a pair of complex
conjugate points in $\projc \smin \projr$ up to the action of $PSL_{2}(\R)$.
Due to Lemma \ref{lem_conf_sp},
 $\Gs_2$-invariant curves $\curve \in \real{M_{\BS}^{\Gs_2}}$  are
 identified with tuples $(\im,-\im, x_{s_3}, \infty), -\infty < x_3 <\infty$.
Hence, $\real{M_{\BS}^{\Gs_2}}= \projr \smin \{\infty\}$ and its
compactification $\real{\konj{M}_{\BS}^{\Gs_2}}$ is $\projr$. 

The main stratum $\real{M_{\BS}^{\Gs_3}}$ of $\real{\konj{M}_{\BS}^{\Gs_3}}$
has different pieces parameterizing $\Gs_3$-invariant   curves
with non-empty and empty  real parts: The
elements of the subspace of $\real{M_{\BS}^{\Gs_3}}$ parameterizing
$\Gs_3$-invariant  curves with $\real{\GS} \not= \emp$ are identified
with tuples $(\Gl \im,\im, -\Gl \im, -\im)$ where $\Gl \in ]-1,1[ \smin \{0\}$.
Similarly, the elements of the subspace of $\real{M_{\BS}^{\Gs_3}}$
parameterizing $\Gs_3$-invariant  curves with $\real{\GS}= \emp$ are  identified
with tuples  $(\Gl \im, \im, - \im/ \Gl,- \im)$, where
$\Gl \in ]-1,1[$. Note that, the strata parameterizing $\real{\GS}
\not= \emp$ and $\real{\GS} = \emp$ are joined through the boundary
points corresponding to $\Gs_3$-invariant curves of type 3.
The compactification of $\real{M_{\BS}^{\Gs_3}}$ is again $\projr$.
\end{ex}

\begin{ex} \label{exa_strata2}
Let $\BS = \{s_1,\cdots,s_5\}$. The moduli space $\cmod{\BS}$ has
three different types of real structures that are given by
\begin{eqnarray*}
\Gs_1 = \mathrm{id}, \
\Gs_2 = \left(
           \begin{array}{ccccc}
              s_1 & s_2 & s_3 & s_4 & s_5 \\
              s_2 & s_1 & s_3 & s_4 & s_5 \\
           \end{array}
         \right),  \
\Gs_3 = \left(
         \begin{array}{ccccc}
              s_1 & s_2 & s_3 & s_4 & s_5\\
              s_3 & s_4 & s_1 & s_2 & s_5\\
         \end{array}
\right).
\end{eqnarray*}
All other real structures of $\cmod{\BS}$ are conjugate to
$c_{\Gs_i}$ since the automorphism group of $\cmod{\BS}$ is $\S_5$
in this case  (see \cite{dik}).

The main stratum $\real{M_{\BS}^{\Gs_1}}$  is identified with the
configuration space of five distinct points in $\projr$ modulo
$PSL_{2}(\R)$, following Lemma \ref{lem_conf_sp}. It is $(\projr \setminus \{0,1,\infty\})^{2}
\setminus \GD$, where $\GD$ is union of all diagonals. Each
connected component of $\real{M_{\BS}^{\Gs_1}}$ is isomorphic to
 a two-dimensional open simplex.
The closure of each cell can be obtained by adding the boundaries
described in Theorem \ref{thm_strata}
(for an example see Fig. \ref{fig_rmoduli3}a).
The cells corresponding to u-planar
trees $\BGt_1$ and $\BGt_2$ are glued along their common stratum
corresponding to $\BGg$, which gives $\BGt_i,i=1,2$ by contracting an edge of
$\BGg$.
The moduli space  $\real{\konj{M}_{\BS}^{\Gs_1}}$ is  a
torus with 3 points blown up (see Fig. \ref{fig_rmoduli3}b).
\begin{figure}[htb]
\centerfig{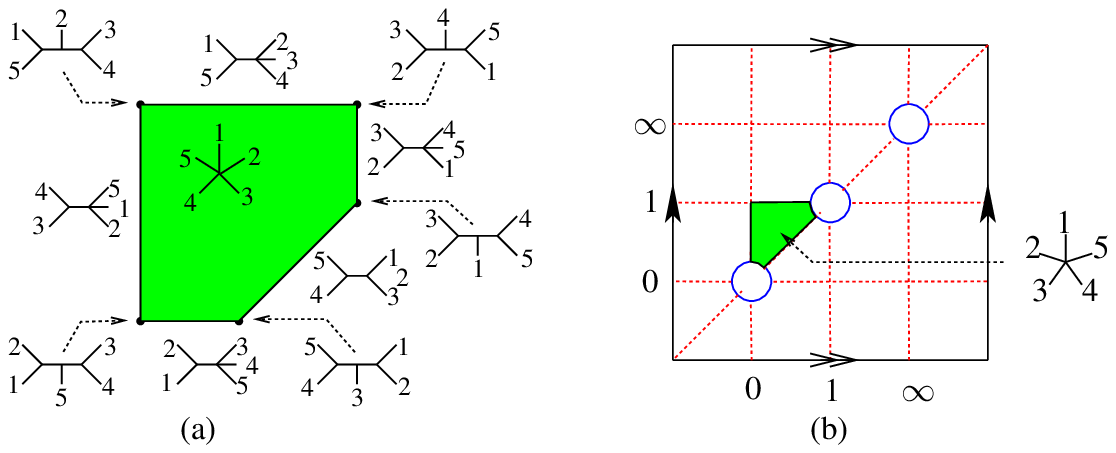,height=1.8 in} %
\caption{(a) Stratification of  $C_{\kBGt}$.
(b) The stratification of $\rmod{\BS}$ for $|\BS|=5$ and $\Gs=\mathrm{id}$.} %
\label{fig_rmoduli3}
\end{figure}

The main stratum $\real{M_{\BS}^{\Gs_2}}$  is diffeomorphic to the
configuration space of a conjugate pair of points on $\projc$.
Following Lemma \ref{lem_conf_sp},  $\Gs_2$-invariant curves in
$\real{M_{\BS}^{\Gs_2}}$ are identified with tuples
$(z,\bar{z},0,1,\infty)$ where $z \in \C \smin \R$; i.e., the main
stratum is $\projc \smin \projr$. The moduli space
$\real{\konj{M}_{\BS}^{\Gs_2}}$
is  a sphere with three points blown up according to the
stratification given in Theorem \ref{thm_strata}.

Finally, $\Gs_3$-invariant curves in $\real{M_{\BS}^{\Gs_3}}$ are identified with
the configurations of a conjugate pair of points in
$\H^+ \smin \{\im\} \sqcup \H^-  \smin \{-\im\}$ following \ref{lem_conf_sp}.
These configurations are given by
$(z,\im,\bar{z},-\im,\infty)$, and the main space $\real{M_{\BS}^{\Gs_3}}$ is
identified with $\projc \smin (\projr \cup \{\im,-\im\})$. The moduli space
 $\real{\konj{M}_{\BS}^{\Gs_3}}$ is a sphere with a point blown up.
\end{ex}


\section{Homology of the strata of $\rmod{\BS}$}
\label{ch_rel_hom}

In this section,  we calculate the homology of the strata of
$\rmod{\BS}$ relative to the union of their substrata of codimension
one and higher.

\subsection{Forgetful morphism revisited}
\label{sec_re_forget}

This section discusses the restriction of the forgetful morphism to
$\rmod{\BS}$.

\subsubsection{The forgetful morphism and u-planar trees}

Let $\BS' \subset \BS$ such that $\Gs(\BS')=\BS'$. Denote the
restriction $\Gs$ on $\BS'$ by $\Gs'$. In this case, the morphism $\pi_{\BS \smin
\BS'}: \cmod{\BS} \to \cmod{\BS'}$ forgetting the points labelled by
$\BS \smin \BS'$ is a real morphism; i.e., $\pi_{\BS \smin \BS'}
\circ c_{\Gs} = c_{\Gs'} \circ \pi_{\BS \smin \BS'}$. Therefore,
$\pi_{\BS \smin \BS'}$ maps the real part of $(\cmod{\BS},c_\Gs)$
onto the real part of $(\cmod{\BS'},c_{\Gs'})$.

Let $\BGg^*$ be the u-planar tree of  $(\GS^*;\ve{p}^*) \in \rmod{S}$, and
let $\BGg$ be the u-planar tree of $\pi_{\BS \smin \BS'}((\GS^*;\ve{p}^*))$.
Then, we say that $\BGg$ is obtained by {\it
forgetting}  the tails $\BS \smin \BS'$ of  $\BGg^*$.


\subsubsection{Forgetting a conjugate pair of labeled points}
\label{sec_unutma}

Let $\BS = \{s_1,\cdots,s_n\}$ and $\Gs \ne \mathrm{id}$. Let $\BS' = \BS
\smin \{s,\bar{s}\}$ for a pair $s, \bar{s} := \Gs(s) \in
\bconj(\Gs)$. Let $\Gs'$ be the restriction of $\Gs$ on $\BS'$. From
now on, we denote the morphism $\pi_{\{s,\bar{s}\}}:
\rmod{\BS} \to \rmods{\BS'}$, forgetting the labeled points
$p_s,p_{\bar{s}}$  by $\pi$.

Let $C_{\BGg^*}$ be a stratum of $\rmod{\BS}$, and identify it with
$C_{\kBGg^*}$ where $\kBGg^*$ is an o-planar representative of $\BGg$.
Let $\pi: C_{\kBGg^*} \to C_{\kBGg}$ be the restriction of the forgetful
map $\pi$ to $C_{\kBGg^*}$. Let $v_s$ be the vertex adjacent to the tail
$s$ i.e., $v_s :=\Bdd_{\Gg^*}(s)$.
Whenever $v_s \in \BV^\R_{\Gg^*}$ and $|v_s|=3$,
there is a unique vertex in $\BV^\R_{\Gg^*}$ next to $v_s$ since
$v_s$ supports both $s$ and $\bar{s}$ and a real edge connecting
$v_s$ to the rest of   $\kBGg^*$. We denote this vertex closest to
$v_s$ by $v_c$.

We will denote the fibers of the forgetful map
$\pi: C_{\kBGg^*} \to C_{\kBGg}$   by $A_{\kBGg^*}$.

\begin{lem}
\label{lem_fil_I}

(a) Let $\curve \in C_{\kBGg}$ be a $\Gs'$-invariant curve of type 1,
and $s \in \BF^+_{\Gg^*}$ (resp. $s \in \BF^-_{\Gg^*}$).  Then, the
fiber  $A_{\kBGg^*}$ over $\curve$ is

\begin{itemize}

\item[1.] a two-dimensional open disc $\GS^+_{v_s}$ (resp. $\GS^-_{v_s}$)
minus the special points $p_f$ where $f \in \BF^+_{\Gg^*}(v_s) \smin
\{s\}$ (resp. $f \in \BF^-_{\Gg^*}(v_s) \smin \{s\}$) if  $v_s \in
\BV_{\Gg^*}^\R$ and $|\BF_{\Gg^*}(v_s)| \geq 5$;

\item[2.] a two-dimensional  sphere $\GS_{v_s}$ minus the special points
$p_f$ where $f \in \BF_{\Gg^*}(v_s) \smin \{s\}$  if  $v_s \not\in
\BV_{\Gg^*}^\R$ and $|\BF_{\Gg^*}(v_s)| \geq 4$;

\item[3.] an open interval if  $v_s \in  \BV_{\Gg^*}^\R$ and
$|\BF_{\Gg^*}(v_s)| = 4$;

\item[4.] an open interval if  $v_s \in  \BV_{\Gg^*}^\R$,
$|\BF_{\Gg^*}(v_s)| = 3$, $|\BF_{\Gg^*}(v_c)| \geq 4$ and
$|\BF^\R_{\Gg^*}(v_c)| > 1$;

\item[5.] a circle if  $v_s \in  \BV_{\Gg^*}^\R$, $|\BF_{\Gg^*}(v_s)| = 3$,
$|\BF_{\Gg^*}(v_c)| \geq 4$ and $|\BF^\R_{\Gg^*}(v_c)| = 1$;

\item[6.] a point if  $v_s \in  \BV_{\Gg^*}^\R$ and
$|\BF_{\Gg^*}(v_s)| = |\BF_{\Gg^*}(v_c)| =3$;

\item[7.] a point if  $v_s \not\in  \BV_{\Gg^*}^\R$ and
$|\BF_{\Gg^*}(v_s)|  =3$.
\end{itemize}

(b) Let $\curve \in C_{\kBGg}$ be a $\Gs'$-invariant curve of type 2.
Then the fiber $A_{\kBGg^*}$ over $\curve$  is

\begin{itemize}
\item[1.] a two-dimensional  sphere $\GS_{v_s}$ minus the special points
$p_f$ where $f \in \BF_{\Gg^*}(v_s) \smin \{s,\bar{s}\}$  if  $v_s
\in  \BV_{\Gg^*}^\R$ and $|\BF_{\Gg^*}(v_s)| \geq 6$;

\item[2.] a two-dimensional  sphere $\GS_{v_s}$ minus the special points
$p_f$ where $f \in \BF_{\Gg^*}(v_s) \smin \{s\}$  if  $v_s \not\in
\BV_{\Gg^*}^\R$ and $|\BF_{\Gg^*}(v_s)| \geq 4$;

\item[3.] a point if  $|\BF_{\Gg^*}(v_s)|  =3$.
\end{itemize}

(c) Let $\curve \in C_{\kBGg}$ be a $\Gs'$-invariant curve of type 3.
Then the fiber $A_{\kBGg^*}$ over $\curve$  is

\begin{itemize}
\item[1.] a two-dimensional  sphere $\GS_{v_s}$ minus the special points
$p_f$ where $f \in \BF_{\Gg^*}(v_s) \smin \{s\}$  if
$|\BF_{\Gg^*}(v_s)| \geq 4$;

\item[2.] a point if  $|\BF_{\Gg^*}(v_s)|  =3$.
\end{itemize}
\end{lem}

\begin{proof}
Here, we prove only (a). The proofs of all other cases are
essentially the same.

Let $\kBGg$ and $\kBGg^*$ be of type 1. Pick a $\Gs'$-invariant
curve $\curve  \in C_{\kBGg}$. Let $(\GS^*,\ve{p}^*)$ denote the points of
the fiber $A_{\kBGg^*}$ over $\curve$. If $(\GS^*,\ve{p}^*) \in
A_{\kBGg^*}$ does not require a contraction of its component
$\GS^*_{v_s}$ after removing the labeled points $p_s,p_{\bar{s}}$,
then there are two possible subcases due to the stability condition:
\begin{itemize}
 \item[1.] $\GS^*_{v_s}$ is a real component and supports five or
 more special points, or

 \item[2.] $\GS^*_{v_s}$ is not a real component and supports four
 or more special points.
\end{itemize}

If $(\GS^*,\ve{p}^*) \in A_{\kBGg^*}$  requires a contraction
 after removing the labeled points
$p_s,p_{\bar{s}}$, then the component $\GS^*_{v_s}$ supports only
three or four special points. Here, there are three possible
subcases:
\begin{itemize}
 \item[3.] $\GS^*_{v_s}$ is a real component and supports four
 special points,

 \item[4-5-6.] $\GS^*_{v_s}$ is  a real component and supports
 three special points, or

 \item[7.] $\GS^*_{v_s}$ is not a real component and supports
 three special points.
\end{itemize}
In these three cases which require stabilization, we implicitly assumed
that $|\BV_{\Gg^*}|>1$
so as not to violate the stability condition for $\BGg$.

We consider these subcases separately.

We first note that  the positions of the special points lying in the components
other than $\GS_{v_s}$ and $\GS_{v_c}$ are fixed, since we consider
the fiber over a fixed $\Gs'$-invariant curve $\curve$.

1. If  $v_s \in  \BV_{\Gg^*}^\R$ and $|v_s| \geq 5$, then
$(\GS^*,\ve{p}^*) \in A_{\kBGg^*}$ do not require the contraction of
$\GS^*_{v_s}$ after removing $p_s,p_{\bar{s}}$. For all such
$(\GS^*,\ve{p}^*)$, $\GS^*=\GS$, and so $\GS^*_{v_s}=\GS_{v_s}$. Assume
that $s \in \BF^+_{\Gg^*}$. For $f \in \BF^+_{\Gg^*}(v_s) \smin
\{s\}$, the special points $p_f$ are fixed and distinct in $\GS^+_{v_s}$.
Therefore, the elements  $(\GS^*,\ve{p}^*)$ of $ A_{\kBGg^*}$
are determined by the position of the point $p_s$ in $\GS^+_{v_s}$.
Since all special points are distinct, $p_s$ must be in $\GS^+_{v_s}
\smin \bigcup \{p_f\}$; i.e., the fiber is  $\GS^+_{v_s} \smin
\bigcup \{p_f\}$ where $f \in \BF^+_{\Gg^*}(v_s) \smin \{s\}$.

2. If  $v_s \not\in  \BV_{\Gg^*}^\R$ and $|v_s| \geq 4$, then
$(\GS^*,\ve{p}^*) \in  A_{\kBGg^*}$ do not require the contraction of
$\GS^*_{v_s}$ after forgetting $p_s,p_{\bar{s}}$. Similarly to the
above case,  $\GS^*_{v_s}=\GS_{v_s}$, and  $(\GS^*,\ve{p}^*) \in
 A_{\kBGg^*}$ are determined by the position of the point $p_s$ in
$\GS_{v_s}$. Hence, the fiber is  $\GS_{v_s} \smin \bigcup \{p_f\}$
where $f \in \BF^+_{\Gg^*}(v_s) \smin \{s\}$.

3. Since all special points other than $p_s,p_{\bar{s}}$ are fixed, a fiber
of $\pi$ is a family of $\Gs$-invariant curves which, in this case,
is the deformation of the irreducible component
$(\GS^*_{v_s},\ve{p}^*_{v_s})$ with two real special points and the
conjugate pair $p_s,p_{\bar{s}}$. It clearly gives an open interval
(see examples in  Section  \ref{sec_ou_p}).

4-5-6. In these cases, $(\GS^*_{v_s},\ve{p}^*_{v_s})$  cannot be
deformed since $|v_s|=3$. Here, the family of $\Gs$-invariant curves
along the fiber $ A_{\kBGg^*}$ is the deformation of
$(\GS^*_{v_c},\ve{p}^*_{v_c})$ (instead of
$(\GS^*_{v_s},\ve{p}^*_{v_s})$). The fiber parameterizes the nodal
points $\GS^*_{v_s} \cap \GS^*_{v_c}$ which disappear after
forgetting $p_s,p_{\bar{s}}$ and contracting $\GS^*_{v_s}$. Here,
there are three subcases:
(6) The fiber $ A_{\kBGg^*}$ is a point when
$|v_c|=3$ since $(\GS^*_{v_c},\ve{p}^*_{v_c})$ cannot deformed in
this case.
(5) The fiber $ A_{\kBGg^*}$ is a circle when
$\BF^\R_{\Gg^*}(v_c)=1$.  Each $(\GS^*,\ve{p}^*)$ is given
by a different position of the nodal
point $\GS^*_{v_s} \cap \GS^*_{v_c}$ in $\real{\GS_{v_c}}$.
(4) The fiber $ A_{\kBGg^*}$ is
an open interval when $\BF^\R_{\Gg^*}(v_c) > 1$.
Each $(\GS^*,\ve{p}^*)$ is given by a different
position of the nodal point $(\GS^*_{v_s} \cap \GS^*_{v_c}) \in
\real{\GS_{v_c}}$ which can
vary between two other special points lying in the real part of
$\GS^*_{v_c}$.

7. The element $(\GS^*_{v_s},\ve{p}^*_{v_s})$ is unique when $v_s
\not\in  \BV_{\Gg^*}^\R$ and  $|v_s|  =3$, since the contracted
component supports only three special points.

The o-planar trees associated to $(\GS^*,\ve{p}^*)$ are
simply obtained by considering the cases above.
\end{proof}

Consider the forgetful map for closed strata $\pi: \csqc{\kBGg^*}
\to \csqc{\kBGg}$.  We denote the
fiber $\pi^{-1} \curve$ over $\curve \in C_{\kBGg}$ by
$\konj{A}_{\kBGg^*}$ since it is the closure of the fiber $ A_{\kBGg^*}$
of $\pi: C_{\kBGg^*} \to C_{\kBGg}$ over the  same point $\curve$. By using the
stratification of $\csqc{\kBGg^*}$, we obtain a stratification of
the fiber $\konj{A}_{\kBGg^*}$.

\begin{lem}
\label{lem_fil_2}
Let $\konj{A}_{\kBGg^*_i}$ be the fibers of $\pi:
\csqc{\kBGg^*_i} \to \csqc{\kBGg}$ over the same point $\curve \in
C_{\kBGg}$. Then, $\konj{A}_{\kBGg^*_1} \subset \konj{A}_{\kBGg^*_2}$ if
and only if $\kBGg^*_1$ produces $\kBGg^*_2$ by contracting one of
its invariant edges or a pair of conjugate edges.
\end{lem}

\begin{proof}
This statement is a direct corollary of Theorem \ref{thm_strata}.
\end{proof}

\subsection{Homology of the fibers of the forgetful morphism}
\label{sec_f_hom}

Let $\kBGg^*$ be a one-vertex o-planar tree, and let $\pi:
C_{\kBGg^*} \to C_{\kBGg}$ be the morphism  forgetting the labeled points
$p_s,p_{\bar{s}}$, which is discussed in Section \ref{sec_unutma}.
Assume that the fibers are two-dimensional; i.e., a punctured disc or
a punctured sphere (see corresponding cases of Lemma
\ref{lem_fil_I}).

\subsubsection{Case of type 1}

Let $\kBGg^*$ be a one-vertex o-planar tree of type 1. Assume that
$s \in \BF^+_{\Gg^*}$ (resp. $s \in \BF^-_{\Gg^*}$). Then, each
fiber $ A_{\kBGg^*}$ of $\pi$ is homotopy equivalent to a bouquet of
$|\BF^+_{\Gg^*}|-1$ circles $S^1 \vee \cdots \vee S^1$. The
cohomology of $ A_{\kBGg^*}$ is generated by the logarithmic
differentials;
\begin{eqnarray*}
 H^0( A_{\kBGg^*}) &=& \Z, \\
 H^1( A_{\kBGg^*}) &=& \bigoplus_{f } \Z \ \Go_{sf}
\end{eqnarray*}
where
\begin{eqnarray*}
\Go_{sf} = \frac{1}{2 \pi \sqrt{-1}} d \log (z_s - z_f)
\end{eqnarray*}
for $f \in \BF^+_{\Gg^*} \smin \{s\}$ (resp.  $f \in \BF^-_{\Gg^*}
\smin \{s\}$).

The homology with closed support $H_1^c( A_{\kBGg^*})$ is isomorphic to
the cohomology group $H^1( A_{\kBGg^*})$  and generated by the duals of
the logarithmic forms; i.e., by the arcs connecting the punctures
$z_f$ to a point in the boundary of the closure $\konj{A}_{\kBGg^*}$
of the fiber $ A_{\kBGg^*}$.

We denote the dual of the generator $\Go_{sf}$ by $\sE_{sf}$. The
homology group $H_2^c( A_{\kBGg^*})$ is isomorphic to $H^0( A_{\kBGg^*})$
and generated by the relative fundamental class of
$\konj{A}_{\kBGg^*}$. Hence,
\begin{eqnarray*}
 H_2^c( A_{\kBGg^*}) &=& \Z\ [\konj{A}_{\kBGg^*}], \\
 H_1^c( A_{\kBGg^*}) &=& \bigoplus_{f }\Z \ \sE_{sf}
\end{eqnarray*}
where $f \in \BF^+_{\Gg^*} \smin \{s\}$ (resp.  $f \in \BF^{-}_{\Gg^*} \smin \{s\}$).

\subsubsection{Case of type 2}

Let $\kBGg^*$ be a one-vertex o-planar tree of type 2. Then, each
fiber $A_{\kBGg^*}$ of $\pi$ is homotopy equivalent to a bouquet of
$|\BF_{\Gg^*}|-3$ circles $S^1 \vee \cdots \vee S^1$. The
cohomology of $ A_{\kBGg^*}$ is generated by the logarithmic
differentials;
\begin{eqnarray*}
 H^0( A_{\kBGg^*}) &=& \Z, \\
 H^1( A_{\kBGg^*}) &=& \bigoplus_{f } \Z \ \Go_{sf}
\end{eqnarray*}
where
\begin{eqnarray*}
\Go_{sf} = \frac{1}{2 \pi \sqrt{-1}} d \log (z_s - z_f)
\end{eqnarray*}
for $f \in \BF_{\Gg^*} \smin \{s,\bar{s},s_{n}\}$.

The homology with closed support $H_1^c( A_{\kBGg^*})$ is isomorphic to
the cohomology group $H^1( A_{\kBGg^*})$  and generated by the arcs
connecting the pairs of punctures $z_{f_1}, z_{f_2}$. These arcs are
the duals of the cohomology classes $\Go_{sf_1} - \Go_{sf_2}$. We
denote them by $\sF_{s,f_1 f_2}$.  Hence,
\begin{eqnarray*}
 H_2^c( A_{\kBGg^*}) &=& \Z\ [\konj{A}_{\kBGg^*}], \\
 H_1^c( A_{\kBGg^*}) &=& \left( \bigoplus_{f_1 \ne f_2} \Z\ \sF_{s,f_1 f_2} \right) / \cJ_s
\end{eqnarray*}
where $f_i \in \BF_{\Gg^*} \smin \{s,\bar{s}\}$, and the subgroup  $\cJ_s$ is generated by
\begin{eqnarray*}
\sF_{s,f_1 f_2} +\sF_{s,f_2 f_3} + \sF_{s,f_3 f_1}.
\end{eqnarray*}

\subsubsection{Homology of the fibers of $\pi_{\{s\}}: \cmodo{\BS} \to \cmodo{\BS'}$}

Let $|\BS|=n \geq 4$ and $s  \in \BS$
be different than $s_n$.  Let $\BS' = \BS \smin \{s\}$. Then, each
fiber $A_s$ of $\pi_{\{s\}}$ is homotopy equivalent to a bouquet of
$|\BS|-2$ circles $S^1 \vee \cdots \vee S^1$. The cohomology of
$A_s$ is generated by the logarithmic differentials;
\begin{eqnarray*}
 H^0(A_s) &=& \Z, \\
 H^1(A_s) &=& \bigoplus_{f }\Z\ \Go_{sf}
\end{eqnarray*}
where
\begin{eqnarray*}
\Go_{sf} = \frac{1}{2 \pi \sqrt{-1}} d \log (z_s - z_f)
\end{eqnarray*}
for $f \in \BS\smin \{s,s_{n}\}$.

The homology with closed support $H_1^c(A_s)$ is isomorphic to the
cohomology group $H^1(A_s)$  and generated by the arcs connecting
the pairs of punctures $z_{f_1}, z_{f_2}$. These arcs are the duals
of the cohomology classes $\Go_{sf_1} - \Go_{sf_2}$. We denote them
by $\sG_{s,f_1 f_2}$.  Hence,
\begin{eqnarray*}
 H_2^c(A_s) &=& \Z\ [\konj{A}_s], \\
 H_1^c(A_s) &=& \left( \bigoplus_{f_1 \ne f_2} \Z\ \sG_{s,f_1 f_2} \right) / \cJ_s
\end{eqnarray*}
where $f_i \in \BS \smin \{s\}$, and the subgroup  $\cJ_s$ is generated by
\begin{eqnarray*}
\sG_{s,f_1 f_2} +\sG_{s,f_2 f_3} + \sG_{s,f_3 f_1}.
\end{eqnarray*}

\subsection{Homology of the strata}
\label{sec_rel_hom}

\begin{lem}
\label{lem_closed_hom} %
Let $\pi: C_{\kBGg^*} \to C_{\kBGg}$ be the
fibration discussed in Section  \ref{sec_unutma}. Then,
\begin{eqnarray*}
H^c_d (C_{\kBGg^*};\Z) = \bigoplus_{p+q=d} H^c_p (C_{\kBGg};\Z)
\otimes H^c_q (A_{\kBGg^*};\Z).
\end{eqnarray*}
\end{lem}

\begin{proof}
We first consider the subcases where $\dim A_{\kBGg^*} =2$. Assume
that $\kBGg^*$ is of type 1. The strata $C_{\kBGg^*}$ and $C_{\kBGg}$
are given by the products
\begin{eqnarray*}
\prod_{v \in \BV^\R_{\Gg^*}} C_{\kBGg^*_v} \times \prod_{v \in \BV^+_{\Gg^*}} \cmodo{\BF_{\Gg^*}(v)},
\ \ \
\prod_{v \in \BV^\R_{\Gg}}   C_{\kBGg_v}   \times \prod_{v \in \BV^+_{\Gg}}   \cmodo{\BF_{\Gg}(v)}
\end{eqnarray*}
respectively (see, Section \ref{sec_ou_p}). The forgetful map $\pi$
preserves the component $(\GS^*_v,\ve{p}^*_v)$ of $(\GS^*,\ve{p}^*)$
for $v \ne v_s$. Hence, it is the identity map on the factors
\begin{eqnarray*}
C_{\kBGg^*_v}           \to C_{\kBGg_v},
\ \ \
\cmodo{\BF_{\Gg^*}(v)} \to \cmodo{\BF_{\Gg}(v)}
\end{eqnarray*}
for $v \ne v_s$. On the other hand, it gives a fibration
\begin{eqnarray}
\label{eqn_red_fib}
\begin{array}{llcll}
\pi_{res}: & C_{\kBGg^*_{v_s}}         &\to& C_{\kBGg_{v_s}},
& \mathrm{when}\ v_s \in \BV_{\Gg^*}^\R,\ \mathrm{and}\\
\pi_{res}: & \cmodo{\BF_{\Gg^*}(v_s)} &\to& \cmodo{\BF_{\Gg}(v_s)},
& \mathrm{when}\ v_s \not\in \BV_{\Gg^*}^\R,
\end{array}
\end{eqnarray}
with the same fibers $A_{\kBGg^*}$ of $\pi: C_{\kBGg^*} \to C_{\kBGg}$.
Therefore, we only need to consider the fibrations in
(\ref{eqn_red_fib}) to calculate the homology.

The strata $C_{\kBGg^*_{v_s}}$ and $C_{\kBGg_{v_s}}$ are diffeomorphic
to the products of simplices with the products of the upper half plane
minus the diagonals (see, Lemma \ref{lem_conf_sp}).
The map $\pi_{res}$ forgets the coordinate subspace $\H^+$
corresponding to the labeled points $p_s$. For instance, when
$|\BF_{\Gg^*}^\R(v_s)| \geq 3$, the map $\pi_{res}: C_{\kBGg^*_{v_s}} \to
C_{\kBGg_{v_s}}$ is
\begin{eqnarray*}
((\H^+)^{|\BF_{\Gg^*}^+(v_s)|} \smin \GD^*) \times \R^{|\BF_{\Gg^*}^\R(v_s)|-3} \to
((\H^+)^{|\BF_{\Gg}^+(v_s)|} \smin \GD)   \times \R^{|\BF_{\Gg}^\R(v_s)|-3}
\end{eqnarray*}
forgetting the coordinate subspace $\H^+$ of the labeled point
$p_s$. Similarly, $\pi_{res}: \cmodo{\BF_{\Gg^*}(v_s)} \to
\cmodo{\BF_{\Gg}(v_s)}$ is
\begin{eqnarray*}
(\C \smin \{0,1\})^{|\BF_{\Gg^*}^+(v_s)|-3} \smin \GD^*)  \to
(\C \smin \{0,1\})^{|\BF_{\Gg}^+(v_s)|-3} \smin \GD
\end{eqnarray*}
forgetting the coordinate subspace $\C$ of the labeled point $p_s$.

The logarithmic forms $d \log (z_s-z_f)$ give a set of
cohomology classes of the total space $C_{\kBGg^*_{v_s}}$ (resp.
$\cmodo{\BF_{\Gg^*}(v_s)} $). On the other hand, the restrictions of
these logarithmic forms to each fiber generate the cohomology of
that fiber (see, Section \ref{sec_f_hom}). By using the Leray-Hirsch
theorem, we obtain
\begin{eqnarray*}
H^d(C_{\kBGg^*}) = \bigoplus_{p+q=d} H^p(C_{\kBGg}) \otimes
H^q(A_{\kBGg^*}).
\end{eqnarray*}

If $\dim A_{\kBGg^*}= 1$ and the fiber is an open interval, then we
directly have
\begin{eqnarray*}
H^p(C_{\kBGg^*}) = H^p(C_{\kBGg};H^0(A_{\kBGg^*})) = H^p(C_{\kBGg})
\otimes H^0(A_{\kBGg^*}).
\end{eqnarray*}

If $\dim A_{\kBGg^*}= 1$ and the fiber is a circle, then $C_{\kBGg^*}$
is $C_{\kBGg} \times S^1$, and the claim follows from the Kunneth
formula.

If $\dim A_{\kBGg^*}= 0$, then each fiber is a single point and the
statement is obvious.

Finally, the duality between cohomology and homology with closed
support gives us the isomorphisms which we need to complete the
proof.

The same arguments apply to o-planar trees of type 2 and type 3.
\end{proof}

Let $Q_{\BGg}$ be the union of the substrata of $\csqc{\BGg}$ of
codimension one and higher.

\begin{prop}
(a) Let $\BGg$ be an u-planar tree of type 1. The relative homology
group $H_{\dim (C_{\BGg})-d}(\csqc{\BGg},Q_{\BGg})$ is generated by
\begin{eqnarray*}
\sE_{s_{i_1} s_{j_1}}  \otimes \cdots \otimes \sE_{s_{i_d} s_{j_d}}
\end{eqnarray*}
where $j_* < i_*$ and $i_1 < \cdots < i_d \leq |\BF^+_{\Gg}|$. In particular,
\begin{eqnarray*}
H_{\dim (C_{\BGg})}(\csqc{\BGg},Q_{\BGg}; \Z) = \Z \cls{\BGg}
\end{eqnarray*}
where $\cls{\BGg}$ is the relative fundamental class of $\csqc{\BGg}$.

(b) Let $\BGg$ be an u-planar tree of type 2. Then, $H_{\dim
(C_{\BGg})-d}(\csqc{\BGg},Q_{\BGg})$ is generated by
\begin{eqnarray*}
\sF_{s_{i_1} s_{j_1} s_{k_1}}  \otimes \cdots \otimes \sF_{s_{i_d} s_{j_d} s_{k_d}}
\end{eqnarray*}
where $\Gs(s_{j_*}), \Gs(s_{k_*}) \ne s_{i_*}$,  $j_*,k_* < i_*$ and $2< i_1 < \cdots < i_d \leq |\BF_{\Gg}|$. In particular,
\begin{eqnarray*}
H_{\dim (C_{\BGg})}(\csqc{\BGg},Q_{\BGg}; \Z) = \Z \cls{\BGg}
\end{eqnarray*}
where $\cls{\BGg}$ is the relative fundamental class of $\csqc{\BGg}$.

(c) Let $W_\BS$ be the union of the  substrata of $\cmod{\BS}$ of
codimension one and higher. The relative homology group $H_{\dim
(\cmod{\BS})-d}(\cmod{\BS}, W_\BS)$ is generated by
\begin{eqnarray*}
\sG_{s_{i_1} s_{j_1} s_{k_1}}  \otimes \cdots \otimes \sG_{s_{i_d} s_{j_d} s_{k_d}}
\end{eqnarray*}
where  $j_*,k_* < i_*$ and $i_1 < \cdots < i_d \leq |\BS|$. In particular,
\begin{eqnarray*}
H_{\dim (\cmod{\BS}}(\cmod{\BS}, W_\BS; \Z) = \Z [\cmod{\BS}]
\end{eqnarray*}
where $[\cmod{\BS}]$ is the  fundamental class of $\cmod{\BS}$.
\end{prop}

\begin{proof}
We obtain the result by applying the forgetful morphism successively
and using Lemma \ref{lem_closed_hom} and the generators of  the homologies of closed 
support of the fibers given in Section \ref{sec_f_hom}.

It is clear that the top dimensional relative homologies are
generated by the (relative) fundamental classes.
\end{proof}

\section{Graph homology of $\rmod{\BS}$}
\label{ch_graph}

In this section, we give a combinatorial complex whose homology
is the homology of $\rmod{\BS}$.

\subsection{A graph complex of $\rmod{\BS}$: $\bfix (\Gs) \ne \emp$ case}

Let $\Gs \in \S_n$ be an involution such that  $\bfix (\Gs) \ne \emp$.
We define a graded group
\begin{eqnarray}
\cG_d &=& \left( \bigoplus_{\BGg : |\BE_{\Gg}|=|\BS|-d-3}
 H_{\dim (C_{\BGg})} (\csqc{\BGg}, Q_{\BGg};\Z) \right) / I_d, \\
&=& \left( \bigoplus_{\BGg : |\BE_{\Gg}|=|\BS|-d-3} \Z\ \cls{\BGg}
\right) / I_d
\end{eqnarray}
where $\cls{\BGg}$ are  the (relative) fundamental classes of
the strata $\csqc{\BGg}$ of $\rmod{\BS}$.

For $|\bconj(\Gs)| <4$, the subgroup $I_d$ (for degree d) is the trivial
subgroup. In all other cases (i.e., for $|\bconj(\Gs)| \geq 4$),  the subgroup
$I_d$ is generated by the following elements.

\subsubsection*{The generators of the ideal of the graph complex.}

The following subsections $\fR$-1 and $\fR$-2 describe
 the generators of the   ideal of the    graph complex.

\subsubsection*{$\fR$-1. Degeneration of a real vertex.}

Consider an o-planar representative $\kBGg$ of a u-planar
tree $\BGg$ of type 1 such that $|\BE_\Gg|=|\BS|-d-5$, and
consider one of its vertices $v \in \BV_\Gg^\R$ with $|v| \geq 5$
and $|\BF^+_\Gg(v)| \geq 2$. Let $f_i,\bar{f}_i
\in \BF_\Gg \smin \BF_\Gg^\R$ be conjugate pairs of flags for
$i=1,2$ such that $f_1,f_2 \in \BF_\Gg^+(v)$ of $\kBGg$, and
let $f_3 \in \BF_\Gg^\R$. Put $\BF = \BF_\Gg(v) \smin
\{f_1,\bar{f}_1,f_2,\bar{f}_2,f_3\}$

We define two u-planar trees $\BGg_1$ and $\BGg_2$ as follows.

The o-planar representative $\kBGg_1$ of $\BGg_1$
is obtained by inserting a pair of conjugate edges $e=(f_{e},f^{e})$,
$\bar{e} = (f_{\bar{e}},f^{\bar{e}})$ into $\kBGg$ at $v$ in
such a way that $\kBGg_1$ gives $\kBGg$ when we contract the edges $e,\bar{e}$.
Let  $\Bdd_{\Gg_1}(e) =\{\tilde{v},v^{e}\}$, $\Bdd_{\Gg_1}(\bar{e}) =
\{\tilde{v},v^{\bar{e}}\}$. Then, the distribution of flags of $\kBGg_1$
 is given by $\BF_{\Gg_1}(\tilde{v}) = \BF_1
\cup \{f_3, f_e,f_{\bar{e}}\}$, $\BF_{\Gg_1}(v^{e}) = \BF_2
\cup \{f_1,f_2,f^e\}$ and $\BF_{\Gg_1}(v^{\bar{e}}) =
\konj{\BF}_2 \cup \{\bar{f}_1,\bar{f}_2,f^{\bar{e}}\}$
where $(\BF_1, \BF_2, \konj{\BF}_2)$ is an equivariant partition of $\BF$.

The o-planar representative $\kBGg_2$ of  $\BGg_2$ is
obtained by inserting a pair of real edges $e_1=(f_{e_1},f^{e_1})$,
$e_2 =(f_{e_2},f^{e_2})$ into $\kBGg$ at $v$ in
such a way that $\kBGg_2$ produces $\kBGg$ when we contract the edges $e_1, e_2$.
Let $\Bdd_{\Gg_2}(e_1) =\{\tilde{v},v^{e_1}\}$,
$\Bdd_{\Gg_2}(e_2)=\{\tilde{v},v^{e_2}\}$.
The sets of flags of  $\kBGg_2$ are $\BF_{\Gg_2}(\tilde{v}) = \ff{\BF}_1 \cup \{f_3 ,
f_{e_1},f_{e_2}\}$, $\BF_{\Gg_2}(v^{e_1}) = \ff{\BF}_2 \cup
\{f_1,\bar{f}_{1},f^{e_1}\}$, $\BF_{\Gg_2}(v^{e_2}) = \ff{\BF}_3
\cup \{f_2,\bar{f}_{2},f^{e_2}\}$
where $(\ff{\BF}_1, \ff{\BF}_2, \ff{\BF}_3)$ is an equivariant partition of $\BF$.


The u-planar trees $\BGg_1, \BGg_2$ are the equivalence classes
represented by $\kBGg_1, \kBGg_2$ given above.

Then, we define
\begin{eqnarray}
\label{eqn_g_relation1}
\cR(\BGg;v,f_1,f_2,f_3) := \sum_{\BGg_1}  \cls{\BGg_1} -
\sum_{\BGg_2}  \cls{\BGg_2},
\end{eqnarray}
where the summation is taken over all possible $\BGg_i,i=1,2$ for a
fixed set of flags $\{f_1,\bar{f}_1,f_2,\bar{f}_2,f_3\}$.

\subsubsection*{$\fR$-2. Degeneration of a conjugate pair of vertices.}

Consider an o-planar representative $\kBGg$ of a u-planar
tree $\BGg$ of type 1 such that $|\BE_\Gg|=|\BS|-d-5$, and
a pair of its conjugate vertices $v,\bar{v} \in
\BV_\Gg \smin \BV_\Gg^\R$ with $|v|=|\bar{v}| \geq 4$.
Let $f_i \in \BF_\Gg(v), i=1,\cdots,4$ and   $\bar{f}_i \in
\BF_\Gg(\bar{v})$  be the flags conjugate  to $f_i, i=1,\cdots,4$.
Put $\BF = \BF_\Gg(v) \smin
\{f_1,\cdots,f_4\}$. Let $(\BF_1,\BF_2)$ be a two-partition of
$\BF$,  and $\konj{\BF}_1,\konj{\BF}_2$ be the sets of flags that are conjugate
to the flags in $\BF_1,\BF_2$ respectively.

We define two u-planar trees $\BGg_1$ and $\BGg_2$ as follows.

The o-planar representative $\kBGg_1$ of  $\BGg_1$
is obtained by inserting a pair of conjugate
edges $e=(f_{e},f^{e})$, $\bar{e} =(f_{\bar{e}},f^{\bar{e}})$ to
$\kBGg$ at $v,\bar{v}$ in such a way that $\kBGg_1$ produces
$\kBGg$ when we contract the edges $e, \bar{e}$.
Let $\Bdd_{\Gg_1}(e) = \{v_{e},v^{e}\}$,
$\Bdd_{\Gg_1} (\bar{e}) =\{v_{\bar{e}},v^{\bar{e}}\}$.
 The sets of flags of $\kBGg_1$  are $\BF_{\Gg_1}(v_{e}) = \BF_1 \cup \{f_1,f_2,f_e\}$,
$\BF_{\Gg_1}(v^{e}) = \BF_2 \cup \{f_3,f_4,f^e\}$ and
$\BF_{\Gg_1}(v_{\bar{e}}) = \overline{\BF}_1 \cup
\{\bar{f}_1,\bar{f}_2,f_{\bar{e}}\}$, $\BF_{\Gg_1}(v^{\bar{e}}) =
\overline{\BF}_2 \cup \{\bar{f}_3,\bar{f}_4,f^{\bar{e}}\}$.

The o-planar representative $\kBGg_2$ of   $\BGg_2$ is
also obtained by inserting a pair of conjugate edges into $\kBGg$ at the
same vertices $v,\bar{v}$, but the flags are distributed differently on vertices.
Let  $\Bdd_{\Gg_2}(e) =\{v_e,v^e\}$, $\Bdd_{\Gg_2}(\bar{e}) = \{v_{\bar{e}},v^{\bar{e}}\} $.
Then, the distribution of the flags of $\kBGg_2$
 is given by
$\BF_{\Gg_2}(v_{e}) = \BF_1 \cup \{f_1,f_3,f_e\}$,
$\BF_{\Gg_2}(v^{e}) = \BF_2 \cup \{f_2,f_4,f^e\}$,
$\BF_{\Gg_2}(v_{\bar{e}}) = \overline{\BF}_1 \cup \{\bar{f}_1,
\bar{f}_3,f_{\bar{e}}\}$, and
$\BF_{\Gg_2}(v^{\bar{e}}) =
\overline{\BF}_2 \cup  \{ \bar{f}_2,\bar{f}_4,f^{\bar{e}}\}$.

The u-planar trees $\BGg_1, \BGg_2$ are the equivalence classes
represented by $\kBGg_1, \kBGg_2$ given above.

We define
\begin{eqnarray} \label{eqn_g_relation2}
\cR(\BGg;v,f_1,f_2,f_3,f_4) := \sum_{\BGg_1} \cls{\BGg_1}
-  \sum_{\BGg_2}  \cls{\BGg_2},
\end{eqnarray}
where the summation is taken over all $\BGg_1, \BGg_2$ for a fixed
set of flags $\{f_1,\cdots,f_4\}$.

The ideal $I_d$ is generated by $\cR(\BGg;v,f_1,f_2,f_3)$ and
$\cR(\BGg;v,f_1,f_2,f_3,f_4)$ for all $\BGg$ and $v$ satisfying
the required conditions above.

\subsubsection*{The boundary homomorphism of the graph complex}

We define the {\it graph complex} $\cG_\bullet$ of the moduli space
$\rmod{\BS}$ by introducing a boundary map $\dd: \cG_d \to
\cG_{d-1}$
\begin{eqnarray}
\label{eqn_differential} \dd\  \cls{\BGt} \ = \ \sum_{\BGg} \ \pm \
\cls{\BGg},
\end{eqnarray}
where the summation is taken over all u-planar trees $\BGg$ which
give $\BGt$ after contracting one of their real edges.

\begin{thm} \label{thm_homology1}
The homology of the graph complex $\cG_\bullet$ is isomorphic to the
singular homology of $\rmod{\BS}$ for $\bfix(\Gs) \ne \emp$.
\end{thm}

\begin{proof}
First, we note that the statement directly follows when $|\bconj(\Gs)|
= 0$. In this case, the stratification of $\rmod{\BS}$ is a cell
decomposition. Details of this case can be found in \cite{de} and
\cite{ka}. Similarly, the stratifications of $\rmod{\BS}$ are cell
decompositions for the cases  $|\BS|=4$ and $|\bconj(\Gs)|=2$, and  $|\BS|=3$
and $|\bconj(\Gs)|=2$ (see examples in Section \ref{sec_ou_p}).

We prove the statement for the $\Gs \ne \mathrm{id}$ cases by induction on the
cardinality of $\bconj(\Gs)$.

Let $\pi: \rmod{\BS} \to \rmods{\BS'}$ be the morphism forgetting the points
labeled by $s,\bar{s} \in \bconj(\Gs)$. Here, we use the  notations introduced
in Section \ref{sec_unutma}.

Let $B_d$ denote the union of $d$-dimensional strata of
$\rmods{\BS'}$ (i.e, $\bigcup_{\BGg} \csqc{\BGg}$ where
$|\BE_\Gg|=|\BS'|-d-3$). The moduli space $\rmods{\BS'}$ is filtered as
follows
\begin{eqnarray*}
\emp =B_{-1} \subset B_0 \subset \cdots \subset B_{(|\BS'|-4)}
\subset B_{(|\BS'|-3)} = \rmods{\BS'}.
\end{eqnarray*}
The forgetful morphism $\pi$ induces a filtration of $\rmod{\BS}$:
\begin{eqnarray*}
\emp =E_{-1} \subset E_0 \subset \cdots \subset E_{(|\BS'|-4)}
\subset E_{(|\BS'|-3)} = \rmod{\BS}
\end{eqnarray*}
where $E_d = \pi^{-1}(B_d)$.  Then, the spectral sequence of  this
filtration gives us
\begin{eqnarray}
\label{eqn_d-complex} \E^1_{p,q}      = H_{p+q} (E_p,E_{p-1})
\Longrightarrow H_{p+q}(\rmod{\BS};\Z).
\end{eqnarray}
We prove this theorem by writing down the spectral sequence
(\ref{eqn_d-complex})
explicitly. As a first step, we calculate the homology groups
$H_{p+q} (E_p,E_{p-1})$.

From now on, we assume that the statement is true for $\rmods{\BS'}$.

\subsubsection*{Step 1.}

We can write the homology of a pair $(E_p,E_{p-1})$ as a direct sum  of the
homology of its pieces:
\begin{eqnarray*}
H_{p+q} (E_p,E_{p-1}) = \bigoplus_{\BGg : |\BE_\Gg|=|\BS'|-d-3}
H_{p+q} (\pi^{-1}(\csqc{\BGg}),\pi^{-1}(Q_{\BGg})).
\end{eqnarray*}

Consider the following filtration of $\pi^{-1}(\csqc{\BGg})$:
\begin{eqnarray*}
\emp \subset Y_0 \subset Y_1 \subset Y_2 = \pi^{-1}(\csqc{\BGg})
\end{eqnarray*}
where $Y_j$'s are the unions of strata
\begin{eqnarray*}
Y_0 = \bigcup_{\BGg^*_m} \csqc{\BGg^*_m}, \ \ Y_1 =
\bigcup_{\BGz^*_l} \csqc{\BGz^*_l}, \ \ Y_2 = \bigcup_{\BGt^*_k}
\csqc{\BGt^*_k,}
\end{eqnarray*}
such that $\pi$ maps each of these stratum onto $\csqc{\BGg}$, and
the dimension of the fibers of $\pi: Y_j \to \csqc{\BGg}$ is $j$ for
$j=0,1,2$.

By using this filtration, we obtain the following  spectral
sequence
\begin{eqnarray*}
\Y^1_{i,j} = H_{i+j} (Y_i,Y_{i-1} \bigcup_{\BGg} (Y_i \cap
\pi^{-1}(Q_{\BGg}))) \Longrightarrow H_{i+j} (E_p,E_{p-1}).
\end{eqnarray*}
Clearly, $Y_i$ contains strata of dimension $p+i$, and $Y_i \cap
\pi^{-1}(Q_{\BGg})$ contains the substrata that map to $B_{p-1}$
(i.e., substrata of codimension one or higher in $\csqc{\BGg^*}$).
Hence, we have
\begin{eqnarray*}
\Y^1_{0,j} &=&
\bigoplus_{\BGg^*_m}  H_{j}(\csqc{\BGg^*_m},Q_{\BGg^*_m}), \\
\Y^1_{1,j} &=&
\bigoplus_{\BGz^*_l}  H_{j+1}(\csqc{\BGz^*_l},Q_{\BGz^*_l}), \\
\Y^1_{2,j} &=& \bigoplus_{\BGt^*_k}
H_{j+2}(\csqc{\BGt^*_k},Q_{\BGt^*_k}).
\end{eqnarray*}

By using Lemma \ref{lem_closed_hom} (and the isomorphism between
relative homology and homology with closed support), we can write the
groups $\Y_{i,j}^1$ as products of homology groups of the base and
fiber:

The dimension of the fibers $\konj{A}_{\BGg^*_m}$ of $\pi:
\csqc{\BGg^*_m} \to \csqc{\BGg}$ is zero, hence
\begin{eqnarray*}
\Y^1_{0,j} = \bigoplus_{\BGg^*_m} H_j (\csqc{\BGg},Q_{\BGg}) \otimes
H_0^c (A_{\BGg^*_m}).
\end{eqnarray*}

The dimension of the fibers $\konj{A}_{\BGz^*_m}$ of $\pi:
\csqc{\BGz^*_l} \to \csqc{\BGg}$ is one, hence
\begin{eqnarray*}
\Y^1_{1,j} = \bigoplus_{\BGz^*_l}  H_j (\csqc{\BGg},Q_{\BGg})
\otimes H_1^c (A_{\BGz^*_l}).
\end{eqnarray*}

Finally,  the dimension of the fibers $\konj{A}_{\BGt^*_k}$ of $\pi:
\csqc{\BGt_k} \to \csqc{\BGg}$ is two, hence
\begin{eqnarray*}
\Y^1_{2,j} &=&
\bigoplus_{\BGt_k} H_j (\csqc{\BGg},Q_{\BGg}) \otimes H_2^c (A_{\BGt^*_k}), \\
\Y^1_{2,j-1} &=&
\bigoplus_{\BGt_k} H_{j} (\csqc{\BGg},Q_{\BGg})
\otimes H_1^c (A_{\BGt^*_k}),
\end{eqnarray*}
Then, the differentials $d_1^\Y: \Y^1_{2,j} \to \Y^1_{1, j}$ and
$d_1^\Y: \Y^1_{1,j} \to \Y^1_{0, j}$ are given respectively by the
differentials
\begin{eqnarray}
\label{eqn_d11}
\begin{array}{cc}
\dd_*:&  H_2^c  (A_{\BGt^*_k}) \to \bigoplus_{\BGz_l < \BGt_k} H_1^c (A_{\BGz^*_l}), \\
\dd_*:&  H_1^c  (A_{\BGz^*_l}) \to \bigoplus_{\BGg_m < \BGz_l} H_0 (A_{\BGg^*_m}).
\end{array}
\end{eqnarray}
Hence, the differential $d_1^\Y$ maps the fundamental class  $\cls{\BGt^*_k}$ to
$\sum \pm \cls{\BGz^*_l}$, and similarly $\cls{\BGz^*_l}$ to
$\sum \pm \cls{\BGg^*_m}$.

Finally, the differential $d_2^\Y: \Y^1_{2,j} \to \Y_{0, j+1}$ is
given by the differentials
\begin{eqnarray}
\label{eqn_d2} %
\dd_*: H_1^c  (A_{\BGt^*_k}) \to  \bigoplus_{\BGg^*_m < \BGt^*_k}
H_0 (A_{\BGg^*_m}).
\end{eqnarray}
For each pair of zero dimensional fibers $\konj{A}_{\BGg^*_{m_1}}$ and
$\konj{A}_{\BGg^*_{m_2}}$ lying in the same two dimensional fiber
$\konj{A}_{\BGt^*_k}$, there is a
generator in $H_1^c  (A_{\BGt^*_k})$ whose image under $\dd_*$ gives
the difference of these points (see Section \ref{sec_f_hom}).
Therefore, each pair of strata $\csqc{\BGg^*_{m_1}}$,
$\csqc{\BGg^*_{m_2}}$ which  are zero dimensional fibrations over
$\csqc{\BGg}$ are homologous relative to $\pi^{-1} (Q_{\BGg})$ i.e.,
\begin{eqnarray}
\label{eqn_rela_rel} \cls{\BGg^*_{m_1}} - \cls{\BGg^*_{m_2}} =0
\end{eqnarray}
in $H_*(E_p,E_{p-1})$.

It is important to note that  the kernel of the differential
$d_2^\Y$ is trivial. This follows from the fact that the kernel
of $\dd_*$ given in (\ref{eqn_d2}) is trivial, as a consequence of the
homology of the fibers given in Section \ref{sec_f_hom}. Therefore,
the homology $H_*(E_p,E_{p-1})$ is given by the total homology of
the spectral sequence $(\Y_{i,j}/I_0,d_1)$ where $I_o$ is the ideal of
the relations given by (\ref{eqn_rela_rel}).

\subsubsection*{Step 2.}
The calculations in Step 1 imply that the term $\E_{*,*}^1$ is generated
by the relative fundamental classes of the  strata. Moreover, it
admits the relations that are imposed in the definition of
$\cG_\bullet$:

The forgetful morphism maps the chains defined in  $\fR$-1 and
$\fR$-2 for $f_i \ne s$  onto the chains of the same type in $\rmods{\BS'}$.
Hence,  if $f_i \ne s$ for all $i$, then $\cR(\BGg;v,f_1,f_2,f_3)$ and
$\cR(\BGg;v,f_1,\cdots,f_4)$ are homologous to zero since we have assumed
that the statement is true for $\rmods{\BS'}$.

On the other hand, the same statement for $f_1=s$ comes as a consequence of
the calculation of Step 1. For each relation (\ref{eqn_rela_rel}) in
relative homology $H_{*}(E_p,E_{p-1})$, there is a relation in
$H_{*}(\rmod{\BS})$. In fact, the sums (\ref{eqn_g_relation1})
(and also (\ref{eqn_g_relation2}))  are mapped onto
a difference given in
(\ref{eqn_rela_rel}) by the relativization map $rel: E_p \to E_p/E_{p-1}$.

We need to confirm that the sums defined in $\fR$-1 and $\fR$-2 are
indeed homologous to zero. We can show this by using certain
forgetful maps.

{\it The chains of type $\fR$-1.} Consider the
composition of projection $\csqc{\BGt^*} \to \csqc{\BGt_v^*}$ onto
the factor corresponding to vertex $v \in \BV_{\Gt^*}^\R$, and
forgetful map $\csqc{\BGt_v^*} \to \csqc{\BGm_v}$ where $\BGm_v^*$
is a one-vertex o-planar tree with $\BF^+_{\Gm_v} =\{f_1,f_2\}$ and
$\BF^\R_{\Gm_v}=\{f_3\}$ which is obtained by forgetting all tails
of $\BGt_v^*$ but $f_1,\bar{f}_1,f_2,\bar{f}_2,f_3$.

The configuration  space $\csqc{\BGm_v}$ is a two-dimensional disc with a puncture,
and it is stratified as in Figure \ref{fig_five_pts}.
\begin{figure}[htb]
\centerfig{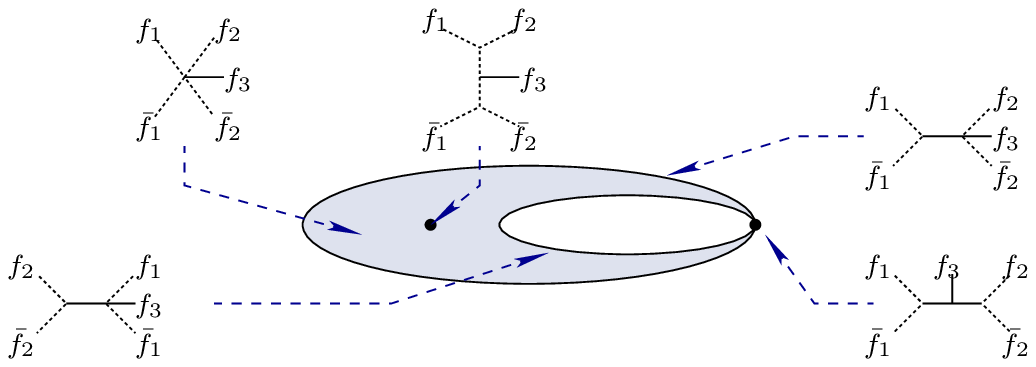,height=1.6in} %
\caption{The strata of $\csqc{\BGm_v}$.} %
\label{fig_five_pts}
\end{figure}

If a codimension two stratum $\csqc{\BGg_1^*}$ of $\csqc{\Gt^*}$ is
in the fiber over the codimension two stratum of $\csqc{\BGm_v}$ which
is lying in its boundary  (see, Figure \ref{fig_five_pts}), then
$\BGg_1^*$ is obtained from $\BGt^*$ by inserting a pair of real
edges at vertex $v$. Similarly, if a codimension two stratum
$\csqc{\BGg_2^*}$ of $\csqc{\Gt^*}$ is in the fiber over the
codimension two stratum of $\csqc{\BGm_v}$ which is lying inside of
$\csqc{\BGm_v}$ (see, Figure \ref{fig_five_pts}), then $\BGg_2^*$ is
obtained from $\BGt^*$ by inserting a pair of conjugate edges at
vertex $v$ as above.  Moreover, since $\csqc{\BGm_v}$ is a punctured
disc, the fibers of the forgetful map over any two points of
$\csqc{\BGm_v}$ are homologous; i.e.,  $\cR(\BGg;v,f_1,f_2,f_3)$ is
homologous to zero.

{\it The chains of type  $\fR$-2.} We show that the sum
(\ref{eqn_g_relation2}) is homologous to zero by using the same
method as in the $\fR$-1 case.
Here, we use a projection map $\csqc{\BGt^*} \to
\cmod{\BF_{\Gt^*}(v)}$. The relations in the complex moduli space
$\cmod{\BF_{\Gt^*}(v)}$ which are introduced by Kontsevich and Manin
induce the relations $\cR(\BGg;v,f_1,\cdots,f_4)=0$ (see, \cite{km2} or \cite{m}).

\subsubsection*{Step 3.} We have a complete description of generators
and relations in $\E^1_{*,*}$. We need to calculate the
differentials.

The first differential  $d_1^\E: \E_{p,q}^1 \to \E_{p-1,q}^1$ is
given by the sum boundary homomorphisms
\begin{eqnarray}
\begin{array}{ccc}
\label{eqn_d1} %
\dd_*: H_{\dim \csqc{\BGt^*}}(\csqc{\BGt^*}) &\to&
\oplus_{\BGg^*} H_{\dim \csqc{\BGg^*}}(\csqc{\BGg^*}) \\
\cls{\BGt^*} &\mapsto & \sum_{\BGg^*} \pm \cls{\BGg^*}
\end{array}
\end{eqnarray}
where $\csqc{\BGt^*} \subset E_p$ and $\csqc{\BGg^*} \subset
E_{p-1}$, and, by contracting a real edge $e$ of each $\BGg^*$, we
obtain $\BGt^*$.

In order to complete the proof, we only need to show that the higher
differentials $d_2^\E$ and $d_3^\E$ of $\E_{*,*}$ vanish.

For dimensional reasons, the differential $d_2^\E$ is zero except
$d_2^\E: \E_{p,1}^1 \to \E_{p-2,2}^1$ and $d_2^\E:\E_{p,0}^1 \to
\E_{p-2,1}^1$.

On the other hand,
\begin{itemize}
\item[{\bf I}] if $\cls{\BGt^*_k} \in \E_{p,2}$,
then either $v_s \in \BV_{\Gt^*_k}^\R$ and $|v_s| \geq 5$,
or  $v_s \not\in \BV_{\Gt^*_k}^\R$ and $|v_s| \geq 4$;

\item[{\bf II}] if $\cls{\BGz^*_l} \in \E_{p,1}$,
then $v_s \in \BV_{\Gz^*_l}^\R$ and either $|v_s|=4$, or
$|v_s|=3$ and $|v_c| \geq 4$;

\item[{\bf III}] if $\cls{\BGg^*_m} \in \E_{p,0}$, then either
$v_s \in \BV_{\Gg^*_m}^\R$ and $|v_s|=|v_c|=3$, or $v_s \not\in
\BV_{\Gg^*_m}^\R$ and $|v_s|=3$
\end{itemize}
(see Lemma \ref{lem_fil_I}).

Now, assume that
\begin{eqnarray*}
d_2^\E (\cls{\BGz^*_l}) = \sum \pm   \cls{\BGt^*_k} \ne 0
\end{eqnarray*}
for $\cls{\BGz^*_l} \in \E_{p,1}^1$. Each $\BGt^*_k$ must
produce $\BGz^*_l$ by contracting one of its real edges, due to
the stratification given in Theorem \ref{thm_strata}.
Note that, the contraction  of  a real edge of $\BGt^*_k$
increases or preserves the valency of the vertex $v_s$.
However, according to condition {\bf II}, the valency $|v_s|$ must
decrease after a contraction  i.e.,  this gives a contradiction.
 Hence, $d_2^\E: \E_{p,1}^1 \to \E_{p-2,2}^1$ must be zero.

By using similar arguments (i.e., comparing valencies of vertices
$v_s$ and $v_c$), we show that  the other higher differentials  $d_2^\E:\E_{p,0}^1 \to
\E_{p-2,1}^1$  and  $d_3^\E:\E_{p,0}^1 \to \E_{p-3,2}^2$
also vanish.

The proof   is completed by observing that the differential
of the graph complex is the sum of the differentials $d_1^\Y$ and
$d_1^\E$ that are defined in (\ref{eqn_d11}),(\ref{eqn_d1})
respectively.
\end{proof}

\begin{rem}
If $|\BS| >4$ and  $|\bfix(\Gs)| \ne 0$, then the moduli space
$\rmod{\BS}$ is not orientable. A combinatorial construction
of the orientation double covering of $\rmod{\BS}$ is given in
\cite{c}. A stratification of the orientation cover is given in
terms of certain equivalence classes of o-planar trees. By following
the same ideas above, it is possible to construct a graph
complex generated by fundamental classes of the strata that
calculates the homology of the orientation double cover of $\rmod{\BS}$.
\end{rem}

\subsection{A graph complex of $\rmod{\BS}$: $\bfix (\Gs) = \emp$ case}

Let $\Gs \in \S_n$ be an involution such that  $\bfix (\Gs) = \emp$.
We define a graded group
\begin{eqnarray}
\sG_d &=& \left( \bigoplus_{\BGg : |\BE_{\Gg}|=|\BS|-d-3}
 H_{\dim (C_{\BGg})} (\csqc{\BGg}, Q_{\BGg};\Z) \right) / I_d, \\
&=& \left( \bigoplus_{\BGg : |\BE_{\Gg}|=|\BS|-d-3} \Z\ \cls{\BGg}
\right) / I_d
\end{eqnarray}
where $\cls{\BGg}$ are   the (relative) fundamental class of
the strata $\csqc{\BGg}$ of $\rmod{\BS}$.

The subgroup $I_d$ (for degree d) is generated by the following
elements.

\subsubsection*{The generators of the ideal of the graph complex.}
\subsubsection*{$\fS$-1. Degeneration of a real vertex: Case of type 1.}

\subparagraph{$\fS$-1.1. $|\BF_\Gg^\R|=0$ case}
Consider an o-planar representative $\kBGg$ of an u-planar
tree $\BGg$ of type 1 such that $|\BE_\Gg|=|\BS|-d-5$ and $|\BF_\Gg^\R|=0$.
Let $v$ be its real vertex, and assume that $|v| \geq 6$. Let 
$f_i \in \BF_\Gg^+, i=1,2,3$,
and let $\bar{f}_i \in \BF_{\Gg}^-$ be their conjugate flags. Put $\BF=
\BF_\Gg(v) \smin \{f_1,f_2,f_3,\bar{f}_1, \bar{f}_2,\bar{f}_3\}$.

We define two u-planar trees $\BGg_1,\BGg_2$ as follows.

The o-planar representative $\kBGg_1$ of  $\BGg_1$
is obtained by inserting  a pair of conjugate
edges $e=(f_{e},f^{e}), \bar{e} = (f_{\bar{e}},f^{\bar{e}})$ into
$\kBGg$ at $v$  in such a way that $\kBGg_1$ produces
$\kBGg$ when we contract the edges $e, \bar{e}$.
Let $\Bdd_{\Gg_1}(e) =\{\tilde{v},v^{e}\}$, $\Bdd_{\Gg_1}(\bar{e})
=\{\tilde{v},v^{\bar{e}}\}$.
The sets of flags are
$\BF_{\Gg_1}(\tilde{v}) = \BF_1 \cup \{f_1,\bar{f}_1,
f_e,f_{\bar{e}}\}$,  $\BF_{\Gg_1}(v^{e}) = \BF_2 \cup
\{f_2,f_3,f^e\}$ and $\BF_{\Gg_1}(v^{\bar{e}}) = \overline{\BF}_2
\cup \{\bar{f}_2,\bar{f}_3,f^{\bar{e}}\}$, where $\BF$ is the
disjoint union of $\BF_1,\BF_2$ and $\overline{\BF}_2$. The set
$\BF_2$ contains the flags that are conjugate to the flags in
$\overline{\BF}_2$ and vice versa.

The o-planar representative $\kBGg_2$ of   $\BGg_2$
is obtained in a similar way. First, we swap
$f_1$  and $\bar{f}_1$ (i.e., put $f_1$ in $\BF_{\Gg}^-$ and
$\bar{f}_1$ in $\BF_{\Gg}^+$). Then, we obtain $\kBGg_2$ by inserting
a pair of conjugate  edges at the vertex $v$ in the same way, but the flags
are distributed differently  $\BF_{\Gg_2}(\tilde{v}) =
\BF_1 \cup \{f_3,\bar{f}_3, f_e,f_{\bar{e}}\}$,
$\BF_{\Gg_2}(v^{e}) = \BF_2 \cup \{\bar{f}_1,f_2,f^e\}$ and
$\BF_{\Gg_2}(v^{\bar{e}}) = \overline{\BF}_2 \cup
\{f_1,\bar{f}_2,f^{\bar{e}}\}$.

The u-planar trees $\BGg_1, \BGg_2$ are the equivalence classes
represented by $\kBGg_1, \kBGg_2$ given above.

Then, we define
\begin{eqnarray}
\label{eqn_son1}
\cR(\BGg;v,f_1,f_2,f_3) := \sum_{\BGg_1}
\cls{\BGg_1} - \sum_{\BGg_2}  \cls{\BGg_2},
\end{eqnarray}
where the summation is taken over all possible $\BGg_i,i=1,2$ for a
fixed set of flags $\{f_i,\bar{f}_i \mid i=1,2,3\}$.

\smallskip
\subparagraph{$\fS$-1.2. $|\BF_\Gg^\R| \ne 0$ case}
The definition of generators of ideal for this case is the same as 
$\fR$-1. Instead of repeating this definition here, we will
refer to (\ref{eqn_g_relation1}) when it is needed.

\subsubsection*{$\fS$-2. Degeneration of a real vertex: Case of type 2.}
Consider an u-planar tree $\BGg$ of type 2 such that
$|\BE_{\Gg}|=|\BS|-d-5$. Let $v$ be its real vertex, and assume that
$|v| \geq 6$. Let $f_i, \bar{f}_i \in \BF_{\Gg}(v)$ be conjugate
pairs of flags for $i=1,2,3$. Put $\BF= \BF_\Gg(v) \smin
\{f_1,f_2,f_3,\bar{f}_1, \bar{f}_2,\bar{f}_3\}$.

We define two u-planar trees $\BGg_1,\BGg_2$ as follows.

The first tree, $\BGg_1$, is obtained by inserting a pair of conjugate
edges $e=(f_{e},f^{e})$, $\bar{e} = (f_{\bar{e}},f^{\bar{e}})$ to
$\BGg$ at $v$ with boundaries $\Bdd_{\Gg_1}(e)
=\{\tilde{v},v^{e}\}$, $\Bdd_{\Gg_1}(\bar{e})
=\{\tilde{v},v^{\bar{e}}\}$. The sets of  flags are given by
$\BF_{\Gg_1}(\tilde{v}) = \BF_1 \cup \{f_1,
\bar{f}_1,f_e,f_{\bar{e}}\}$, $\BF_{\Gg_1}(v^{e}) = \BF_2 \cup
\{f_2,f_3,f^e\}$ and $\BF_{\Gg_1}(v^{\bar{e}}) = \overline{\BF}_2
\cup \{\bar{f}_2,\bar{f}_3,f^{\bar{e}}\}$,  where $\BF$ is a
disjoint union of $\BF_1$, $\BF_2$ and $\overline{\BF}_2$.
$\BF_2$ contains the flags  that are conjugate to the
flags in $\overline{\BF}_2$ and vice versa.

The second tree, $\BGg_2$, is also obtained by inserting a pair of
conjugate edges at the vertex $v$ in the same way, but the flags are
distributed differently on vertices; $\BF_{\Gg_2}(\tilde{v}) = \BF_2
\cup \{f_3, \bar{f}_3,f_e,f_{\bar{e}}\}$, $\BF_{\Gg_2}(v^{e}) =
\BF_1 \cup \{f_1,f_2,f^e\}$ and $\BF_{\Gg_2}(v^{\bar{e}}) =
\overline{\BF}_2 \cup \{\bar{f}_1,\bar{f}_2,f^{\bar{e}}\}$.

Then, we define
\begin{eqnarray}
\label{eqn_son2}
\cR(\BGg;v,f_1,f_2,f_3) := \sum_{\BGg_1}  \cls{\BGg_1} -
\sum_{\BGg_2}  \cls{\BGg_2}.
\end{eqnarray}
Here,  the summation is taken over all possible $\BGg_i,i=1,2$ for a
fixed set of flags $\{f_i,\bar{f}_i \mid i=1,2,3\}$.

\subsubsection*{$\fS$-3. Degeneration of a conjugate pair of vertices.}

Consider an u-planar tree $\BGg$ (of type 1, type 2 or type 3) such
that $|\BE_\Gg|= |\BS|-d-5$, and a pair of its conjugate
vertices $v,\bar{v} \in \BV_\Gg \smin \BV_\Gg^\R$ such that
$|v|=|\bar{v}| \geq 4$. Let $f_i \in \BF_\Gg(v), i=1,\cdots,4$, and
let $\bar{f}_i \in \BF_\Gg(\bar{v})$ be their conjugate flags. Put
$\BF = \BF_\Gg(v) \smin \{f_1,\cdots,f_4\}$. Let $(\BF_1,\BF_2)$ be
a partition of $\BF$,  and $\konj{\BF}_1,\konj{\BF}_2$ be
the sets of conjugate flags.

We define two u-planar trees $\BGg_1$ and $\BGg_2$ as follows.

The first one, $\BGg_1$, is obtained by inserting a pair of conjugate
edges $e=(f_{e},f^{e})$, $\bar{e} =(f_{\bar{e}},f^{\bar{e}})$ into
$\BGg$ at $v,\bar{v}$ such that $\Bdd_{\Gg_1}(e) = \{v_{e},v^{e}\}$,
$\Bdd_{\Gg_1} (\bar{e}) =\{v_{\bar{e}},v^{\bar{e}}\}$. The
sets of flags are given by  $\BF_{\Gg_1}(v_{e}) = \BF_1 \cup \{f_1,f_2,f_e\}$,
$\BF_{\Gg_1}(v^{e}) = \BF_2 \cup \{f_3,f_4,f^e\}$ and
$\BF_{\Gg_1}(v_{\bar{e}}) = \overline{\BF}_1 \cup
\{\bar{f}_1,\bar{f}_2,f_{\bar{e}}\}$, $\BF_{\Gg_1}(v^{\bar{e}}) =
\overline{\BF}_2 \cup \{\bar{f}_3,\bar{f}_4,f^{\bar{e}}\}$.

The second one, $\BGg_2$, is also obtained by inserting a pair of
conjugate edges into $\BGg$ at the same vertices $v,\bar{v}$, but the
flags are distributed differently on vertices. Let  $\Bdd_{\Gg_2}(e)
=\{v_e,v^e\}, \Bdd_{\Gg_2}(\bar{e}) = \{v_{\bar{e}},v^{\bar{e}}\} $. Then,
the sets of flags are given by
$\BF_{\Gg_2}(v_{e}) = \BF_1 \cup \{f_1,f_3,f_e\}$,
$\BF_{\Gg_2}(v^{e}) = \BF_2 \cup \{f_2,f_4,f^e\}$,
$\BF_{\Gg_2}(v_{\bar{e}}) = \overline{\BF}_1 \cup \{\bar{f}_1,
\bar{f}_3,f_{\bar{e}}\}$ and $\BF_{\Gg_2}(v^{\bar{e}}) =
\overline{\BF}_2 \cup  \{ \bar{f}_2,\bar{f}_4,f^{\bar{e}}\}$.

We define
\begin{eqnarray}
\label{eqn_son3}
\cR(\BGg;v,f_1,f_2,f_3,f_4) := \sum_{\BGg_1} \cls{\BGg_1} -
\sum_{\BGg_2}  \cls{\BGg_2},
\end{eqnarray}
where the summation is taken over all such $\BGg_i,i=1,2$ for a fixed
set of flags $\{f_1,\cdots,f_4\}$.

The ideal $I_d$ is generated by the chains $\cR(\BGg;v,f_1,f_2,f_3)$
defined in (\ref{eqn_g_relation1}), (\ref{eqn_son1}) and  (\ref{eqn_son2}),
and $\cR(\BGg;v,f_1,f_2,f_3,f_4)$ defined in  (\ref{eqn_son3}) for all
$\BGg$ and $v$ satisfying the required conditions above.

\subsubsection*{The boundary homomorphism of the graph complex.}

We define the {\it graph complex} $\sG_\bullet$ of the moduli space
$\rmod{\BS}$ by introducing a boundary map $\dd: \sG_d \to
\sG_{d-1}$
\begin{eqnarray*}
\dd\  \cls{\BGt} \ = \ \sum_{\BGg} \ \pm \ \cls{\BGg},
\end{eqnarray*}
where the summation is taken over all u-planar trees $\BGg$ which
give $\BGt$ after contracting one of their real edges.

\begin{thm}
The homology of the graph complex $\sG_\bullet$ is isomorphic to the
singular homology of $\rmod{\BS}$ for $\bfix(\Gs) = \emp$.
\end{thm}

The proof of this theorem is essentially the same as the proof of Theorem
\ref{thm_homology1}. We will not repeat it here.

\section{Fundamental group of $\rmod{\BS}$}
\label{ch_f_group}

In this section, we give a presentation of the fundamental group of
$\rmod{\BS}$ by using the groupoid of paths which are transversal to the
codimension one strata of $\rmod{\BS}$. This idea has been used by Kamnitzer and
Henriques in \cite{hk} to calculate the fundamental group of
$\rmod{\BS}$ for $\Gs=\mathrm{id}$. This section extends their description
to $\pi_1(\rmod{\BS})$ for  an arbitrary involution $\Gs$.

\subsection{Fundamental groups of open parts of strata}

In this section, we consider a particular subset of the set of u-planar
trees. Let $\rtree$ be the set of u-planar trees  having no conjugate pairs
of edges. If $\BGg \in \rtree$ is of type 1, then $\BV_\Gg = \BV_\Gg^\R$.
If $\BGg \in \rtree$ is of type 2, then
$|\BV_\Gg|=1$, and if $\BGg \in \rtree$ is of type 3, then
$|\BV_\Gg|=2$.

For a u-planar tree $\BGg \in \rtree$, the {\it open part}
$\opp{\BGg}$ of the stratum is the closed stratum $\csqc{\BGg}$
minus the union of the closure  of its codimension one strata.

\begin{prop}
\label{cor_homotop} %
For  $\BGg \in \rtree$, the open part of the stratum $\csqc{\BGg}$
is simply connected.
\end{prop}

\begin{proof}
For $\BGg^* \in \rtree$, the open part of a stratum is
\begin{eqnarray*}
\opp{\BGg^*} = \left\lbrace
\begin{array}{ll}
\prod_{v \in \BV_{\Gg^*}} \opp{\BGg^*_v} &
\mathrm{if}\ \BGg^*\ \mathrm{is\ of\ type\ 1}, \\
\opp{\BGg_{v_r}^*} &
\mathrm{if}\ \BGg^*\ \mathrm{is\ of\ type\ 2}, \\
\cmod{\BF_\Gg(v)} &
\mathrm{if}\ \BGg^*\ \mathrm{is\ of\ type\ 3}.
\end{array}
\right.
\end{eqnarray*}
This follows from the fact that the open part $\opp{\BGg^*}$ of
$\csqc{\BGg^*}$ is the product of the open parts of its factors.

Hence, we only need to consider the factors that
correspond to the one-vertex trees.

We prove the statement by induction on the cardinality of
$\bconj(\Gs)$. First, we note that the open parts of the strata of
$\rmod{\BS}$ are contractible for $\BS=\bfix(\Gs)$, $|\BS|=4$ and
$|\bfix(\Gs)|=2$, $|\BS|=4$ and $|\bfix(\Gs)|=0$, and $|\BS|=3$
and $|\bfix(\Gs)|=1$. In these cases, the stratifications are  cell
decompositions, and the open parts of the strata are open discs (see
\cite{de,ka} and examples in Section \ref{sec_ou_p}).

Let $\BGg^*$  be a one-vertex u-planar tree of type 1. Let
$|\BF_{\BGg^*}^+|>0$, and $\pi:\csqc{\BGg^*} \to \csqc{\BGg}$ be the morphism
forgetting the conjugate pairs of points $p_s, p_{\bar{s}}$. Let $O$
be the subset of the fiber $\pi^{-1}\curve$ such that
$(\GS^*,\ve{p}^*) \in O$ does not require any stabilization after
forgetting $p_s, p_{\bar{s}}$. For $(\GS^*,\ve{p}^*) \in O$, $\GS^*=\GS$.
Since all special points are fixed in $\GS^*$, the different elements
of $O$ are given by positions of the labelled point $p_s$. The
labelled point $p_s$ is in $(\GS \smin (\{\mathrm{special\ points}\}
\cup \real{\GS}))/ c_\GS$.  This follows from the fact that all
special points must be distinct (hence, we need to remove special
points and $\real{\GS}$ where $p_s$ and $p_{\bar{s}}$ collide and give a
real node) and $s$ in either $\BF^+_{\Gg^*}$ or $\BF^-_{\Gg^*}$ (so
that, we need to take the quotient with respect to the real
structure $c_{\GS}: \GS \to \GS$).

A degeneration of   $(\GS^*,\ve{p}^*) \in O$, which is
obtained as the limit as $p_s$ goes to a special point in $\GS \smin
\real{\GS}$, gives us an element in $\opp{\BGg^*}$, since the limit
element has an additional  conjugate pair of edges. On the other
hand, a degeneration of  $(\GS^*,\ve{p}^*)$, which is obtained as the
limit as $p_s$ goes to a point in $\real{\GS}$, gives a curve with a real
node; i.e., this limit does not lie in $\opp{\BGg^*}$. Therefore, the
restriction of the forgetful map $\pi: \opp{\BGg^*} \to \opp{\BGg}$
has a fiber $(\GS \smin \real{\GS})/c_\GS$ over $\curve \in
\opp{\BGg}$. It is clear that the fiber is simply connected.

If we assume simply connectedness of $\opp{\BGg}$, then
$\opp{\BGg^*}$ is clearly simply connected. We prove the statement
by induction on the cardinality of the labeling set $\bconj(\Gs)$.

The proofs for u-planar trees $\BGg \in \rtree$ of type 2 and type 3
are the same as for the type 1 case.  The fiber of the forgetful map $\pi:
\opp{\BGg^*} \to \opp{\BGg}$ over $\curve$ is $\GS$ when $\BGg^*$ is
of type 2, and $\GS/c_\GS$ when  $\BGg^*$ is of type 3. In both
cases, the fibers are simply connected.
\end{proof}

\begin{prop}
The moduli space $\rmod{\BS}$ is stratified by simply connected
subspaces $\opp{\BGg}$.
\end{prop}

\begin{proof}
This directly follows from the fact that the open parts $\opp{\BGg}$
of the strata $\csqc{\BGg}$ are pairwise disjoint.
\end{proof}

\subsection{A groupoid of paths in $\rmod{\BS}$}
We consider the following groupoid $\cP$ of paths in $\rmod{\BS}$.
 Choose an element $(\GS(i),\ve{p}(i))$ in every
connected component $C_{\BGt_{(i)}}$ of the main stratum $\rmodo{\BS}$.

The objects Ob($\cP$) are these elements $(\GS(i),\ve{p}(i))$.

The morphisms of $\cP$ are given as follows:
Let $\csqc{\BGt_{(i)}}$ and $\csqc{\BGt_{(j)}}$ be a pair
of top-dimensional adjacent strata in $\rmod{\BS}$, and
let $\csqc{\BGg_{ij}}$ be their common codimension
one stratum. The morphism $\langle\BGg_{ij} \rangle$
is the homotopy class of oriented paths in $\rmod{\BS}$
that start from the point $(\GS(i),\ve{p}(i))$, pass
through the common codimension one (open) stratum $C_{\BGg_{ij}}$, and end at
$(\GS(j),\ve{p}(j))$. Notice that such paths connecting
a pair of  points $(\GS(i),\ve{p}(i))$ and
$(\GS(j),\ve{p}(j))$ are homotopic, since the open parts of $\csqc{\BGt_{(i)}}$,
$\csqc{\BGt_{(j)}}$ and $\csqc{\BGg_{ij}}$ are simply
connected (see Proposition \ref{cor_homotop}). The  homotopy classes
of paths that intersect only codimension 1 strata are given
by concatenations of  $\langle\BGg_{ij} \rangle$'s.

\begin{thm} \label{thm_fund}
The fundamental group $\pi_1(\rmod{\BS})$ is
generated by the loops
\begin{eqnarray}
\label{eqn_loop}%
\langle \BGg_{i_1 i_2}  \ \BGg_{i_2 i_3} \cdots \ \BGg_{i_{n-1} i_n} \ \BGg_{i_n i_1} \rangle
: =
\langle \BGg_{i_1 i_2} \rangle
\langle \BGg_{i_2 i_3} \rangle \cdots
\langle \BGg_{i_n i_1} \rangle
\end{eqnarray}
which are subject to the following relations:

\begin{enumerate}
\item Let $\BGg_{i_1 i_2}$ be a u-planar tree having only one edge, then
\begin{eqnarray}
\label{eqn_g_rel1} %
\langle \BGg_{i_1 i_2}  \ \BGg_{i_2 i_1} \rangle
=1.
\end{eqnarray}

\item Let $\BGg'$ with $|\BE_{\Gg'}|=2$, and let
$\BGg_{ij}, \BGt_{ji}$ be the only u-planar
trees that are obtained by contracting one of the edges of $\BGg'$, then
\begin{eqnarray}
\label{eqn_g_rel2} %
\langle \BGg_{ij}  \ \BGt_{ji} \rangle
=1.
\end{eqnarray}

\item Let $\BGg'$ with $|\BE_{\Gg'}|=2$, and let $\BGg_{i_1 i_2}, \cdots,
\BGg_{i_4 i_1}$ be the u-planar
trees that are obtained by contracting an edge of $\BGg'$, then
\begin{eqnarray}
\label{eqn_g_rel22} %
\langle \BGg_{i_1 i_2}  \ \BGg_{i_2 i_3}  \ \BGg_{i_3 i_4}  \ \BGg_{i_4 i_1}  \rangle
=1.
\end{eqnarray}
\end{enumerate}
\end{thm}

\begin{proof} Every loop in $\rmod{\BS}$ is homotopic to a loop which
is transversal to codimension one strata. Such transversal loops can be obtained by
perturbing the original loops. Hence, we can choose the loops given in
(\ref{eqn_loop}) as representatives of the homotopy classes of loops.

These loops are subject to the following relations:

\begin{itemize}
\item[1.]
The concatenation
of a path $\langle \BGg_{i_1 i_2} \rangle$ with the reverse path
$\langle \BGg_{i_2 i_1} \rangle$ is obviously homotopic to a
point and gives the relation (\ref{eqn_g_rel1}).
\end{itemize}

On the other hand, if two paths in $\rmod{\BS}$ are homotopic, then they are homotopic by a
homotopy of paths that are transversal to the codimension one strata. Therefore, the
homotopy relations arise from the passing of the paths through codimension two strata:
Let $\BGg'$ be a u-planar tree corresponding to a codimension two stratum
of $\rmod{\BS}$. The stratum $\csqc{\BGg'}$ is contained  either in
two or in four codimension one strata, since we can obtain two or four u-planar
trees by contracting one of the two edges of $\BGg'$.

\begin{itemize}
\item[2.]
If only two codimension one strata, $\csqc{\BGg_{ij}}$ and $\csqc{\BGt_{ji}}$,  intersect along
the codimension two stratum $\csqc{\BGg'}$, then  the concatenation of
$\langle \BGg_{ij}  \rangle$ and $\langle \BGt_{ji} \rangle$ gives the loop
\begin{eqnarray*}
\langle \BGg_{ij}  \ \BGt_{ji} \rangle
\end{eqnarray*}
around  $\csqc{\BGg'}$ which is contractible; i.e., gives relations in  (\ref{eqn_g_rel2}).
\end{itemize}

\begin{itemize}
\item[3.]
If there are four codimension one strata
$\csqc{\BGg_{i_1 i_2}}, \cdots, \csqc{\BGg_{i_4 i_1}}$ intersecting
along the codimension two stratum $\csqc{\BGg'}$, then the loop
\begin{eqnarray*}
\langle \BGg_{i_1 i_2}  \ \BGg_{i_2 i_3}  \ \BGg_{i_3 i_4}  \ \BGg_{i_4 i_1}  \rangle
\end{eqnarray*}
around  $\csqc{\BGg'}$ is contractible and gives the relations (\ref{eqn_g_rel22}).
\end{itemize}
\end{proof}

\begin{rem}
The fundamental group of $\rmod{\BS}$ is particularly interesting for $\Gs =\mathrm{id}$,
since  $\rmod{\BS}$ is $K(\pi,1)$ in this case. However, this is not true
when $\Gs \ne \mathrm{id}$. For instance, if $|\bconj(\Gs)| =6$, then the open part of each
top dimensional stratum is topologically a 2-dimensional sphere i.e., higher homotopy
groups of $\rmod{\BS}$ cannot be trivial in this case.
\end{rem}

\begin{rem}
A similar presentation of the fundamental group of the orientation double cover
of $\rmod{\BS}$ is given in \cite{th}.
\end{rem}



\end{document}